\pdfoutput=1
\documentclass[english]{jnsao}

\usepackage[utf8]{inputenc}
\usepackage{preamble}

\usepackage{enumitem}
\setlist[enumerate,1]{label=(\roman*)}

\usepackage{cleveref}

\theoremstyle{definition}
\newtheorem{assumption}[theorem]{Assumption}

\providecommand\given{\nonscript\;\delimsize|\nonscript\;\mathopen{}}
\DeclarePairedDelimiterX\set[1]\{\}{#1}
\DeclarePairedDelimiterX\seq[1](){#1}

\DeclarePairedDelimiterX\dual[2]{\langle}{\rangle}{#1,#2}
\DeclarePairedDelimiterX\innerprod[2](){#1,#2}

\DeclarePairedDelimiter\abs{\lvert}{\rvert}
\DeclarePairedDelimiter\norm{\lVert}{\rVert}

\DeclarePairedDelimiter\parens()
\DeclarePairedDelimiter\bracks[]
\DeclarePairedDelimiter\braces\{\}

\newcommand\N{\mathbb{N}}
\newcommand\R{\mathbb{R}}
\newcommand\barR{\bar\R}

\renewcommand\d{\mathop{}\!\mathrm{d}}
\newcommand\dx{\d x}
\newcommand\dt{\d t}
\newcommand\dxt{\d(x,t)}

\newcommand{\weakly}{\rightharpoonup}
\newcommand{\weaklystar}{\stackrel\star\rightharpoonup}

\newcommand{\dualspace}{^\star}
\newcommand{\bidualspace}{^{\star\star}}

\newcommand{\firstsubderivative}{^\downarrow}
\newcommand{\weakfirstsubderivative}{^{\sim}}

\newcommand\KK{\mathcal{K}}

\newcommand\Uad{U_{\mathrm{ad}}}
\newcommand\tangentcone[1]{\mathcal{T}_{\Uad}\parens*{#1}}
\newcommand\normalcone[1]{\mathcal{N}_{\Uad}\parens*{#1}}
\newcommand\criticalcone[1]{C_{#1}}
\newcommand\LinftyOmegaT{L^\infty(\Omega_T)}
\newcommand\LtwoOmegaT{L^2(\Omega_T)}
\newcommand\indicatorofUad{\delta_{\Uad}}
\newcommand\characteristicfunction[1]{\chi_{#1}}
\newcommand\projection[2]{\mathrm{P}_{#1}\parens*{ #2 }}
\newcommand\Ltwonorm[1]{\norm{#1}_{\LtwoOmegaT}}
\newcommand\normLtwotime[1]{\norm{#1}_{L^2(0, T)}}
\newcommand\normLoneOmega[1]{\norm{#1}_{L^1(\Omega)}}
\newcommand\jone[1]{\norm{#1}_{L^1(\Omega_T)}}
\newcommand\jtwo[1]{\norm{#1}_{L^2(0,T;L^1(\Omega))}}
\newcommand\jthree[1]{\norm{#1}_{L^1(\Omega;L^2(0,T))}}
\DeclarePairedDelimiterXPP\measure[2]{\lambda^{#1}}(){}{#2}

\newcommand\OmegaT{{\Omega_T}}

\DeclareMathAlphabet{\mathpzc}{OT1}{pzc}{m}{it}
\newcommand\oo{\mathpzc{o}}

\newcommand\Sign{\operatorname{Sign}}
\newcommand\dom{\operatorname{dom}}

\manuscriptsubmitted{2023-11-28}
\manuscriptaccepted{2024-11-15}
\manuscriptvolume{5}
\manuscriptnumber{12604}
\manuscriptyear{2024}
\manuscriptdoi{10.46298/jnsao-2024-12604}
\manuscriptlicense{CC-BY-NC-SA 4.0}

\begin{document}
\title{%
	Second-order conditions
	for
	spatio-temporally sparse optimal control
	via
	second subderivatives
}

\author{%
	Nicolas Borchard%
	\thanks{%
		Brandenburgische Technische Universität Cottbus-Senftenberg,
		Institute of Mathematics,
		03046 Cottbus, Germany,
		\email{nicolas.borchard@b-tu.de},
		\orcid{0009-0007-7358-7737}%
	}
	\and
	Gerd Wachsmuth%
	\thanks{%
		Brandenburgische Technische Universität Cottbus-Senftenberg,
		Institute of Mathematics,
		03046 Cottbus, Germany,
		\email{gerd.wachsmuth@b-tu.de},
		\orcid{0000-0002-3098-1503}%
	}
}
\shortauthor{Borchard, Wachsmuth}

\maketitle

\begin{abstract}
	We address
	second-order optimality conditions
	for optimal control problems
	involving sparsity functionals
	which induce spatio-temporal sparsity patterns.
	We employ the notion of (weak) second subderivatives.
	With this approach,
	we are able to reproduce the results from
	Casas, Herzog, and Wachsmuth (ESAIM COCV, 23, 2017, p.\ 263--295).
	Our analysis yields a slight improvement
	of one of these results
	and also opens the door for the sensitivity analysis
	of this class of problems.
\end{abstract}


\begin{msc}
	\mscLink{49K27},
	\mscLink{49K20}
\end{msc}

\section{Introduction}
\label{sec:introduction}
We are interested in second-order optimality condition for optimization problems of the form
\begin{equation}
	\label{eq:prob}
	\text{Minimize}
	\quad
	F(u) + \mu j_i(u)
	\qquad\text{w.r.t.\ } u \in \Uad,
\end{equation}
where $F \colon \Uad \to \R$
is assumed to be smooth,
\begin{equation*}
	\Uad :=
	\set{
		u \in L^2( \Omega \times (0,T) )
		\given
		\alpha \le u \le \beta
		\; \text{a.e.}
	}
\end{equation*}
is the feasible set
with constants $\alpha < \beta$,
$\mu > 0$ is a scaling parameter
and $j_i$ is one of the sparsity functionals
\begingroup
\allowdisplaybreaks
\begin{subequations}
	\label{eq:js}
	\begin{align}
		\label{eq:js:1}
		j_1(u) &:= \jone{u} := \int_{\Omega \times (0,T)} \abs{u(x,t)} \dxt, \\
		j_2(u) &:= \jtwo{u} := \bracks*{ \int_0^T \normLoneOmega{u(\cdot, t)}^2 \dt } ^ {1 / 2}, \\
		j_3(u) &:= \jthree{u} := \int_\Omega \normLtwotime{u(x, \cdot)} \dx
		.
	\end{align}
\end{subequations}
\endgroup
For a special choice of $F$,
problem \eqref{eq:prob}
has been considered in
\cite{CasasHerzogWachsmuth2015:1}.
Therein, the smooth part $F$ involves the solution map of a semilinear parabolic equation,
see \cref{sec:application} for details.

Sufficient optimality conditions of second-order
typically provide quadratic growth in the neighborhood of a minimizer.
This is of uttermost importance for the stability of the minimizer under perturbations,
for the convergence of optimization methods
and for the numerical analysis of discretization schemes,
see \cite{CasasTroeltzsch2014}
and the references therein.
Another important role is played by
necessary conditions of second-order.
They will be satisfied by every local minimizer
and the ``distance'' between
the necessary and the sufficient conditions
can serve as a judgement concerning their sharpness.

There is a rich literature for second-order conditions for PDE-constrained
optimal control problems.
We refer again to \cite{CasasTroeltzsch2014}
and the references therein.
The first contribution concerning second-order conditions
for optimal control problems with nonsmooth sparsity functionals
is \cite{CasasHerzogWachsmuth2012}.
Therein, the $L^1$-norm (as in \eqref{eq:js:1}) was considered,
which has a significantly simpler structure
than the functionals $j_i$ from \eqref{eq:js}.
As already said, the other functionals from \eqref{eq:js}
were already investigated in 
\cite{CasasHerzogWachsmuth2015:1}.
Further contributions which address
more complicated state equations
are
\cite{CasasRyllTroeltzsch2015,SprekelsTroeltzsch2023}.

The structure of second-order conditions
depend crucially on the presence of an $L^2$-Tikhonov term in the objective.
If such a term is included,
the problem can be analyzed in $L^2$,
one obtains a quadratic growth condition in $L^2$ for the control
and the gap between the necessary and sufficient conditions of second order is as small as possible,
see, e.g., \cite{CasasTroeltzsch2012}.
In absence of the Tikhonov term,
the first contribution for PDE control is
\cite{Casas2012}.
Therein,
the author proves a second-order sufficient condition
which guarantees a quadratic growth in the state variable,
but necessary conditions of second order are not available.
Under a structural assumption on the adjoint state,
\cite{CasasWachsmuthWachsmuth2017}
provides a second-order condition which gives quadratic growth in $L^1$
for the control, again, no necessary conditions of second order were addressed.
By using the notion of second subderivatives,
it was demonstrated in \cite{ChristofWachsmuth2017:1},
that (under the structural assumption)
the problem should be analyzed in the space of measures
and the authors were able to characterize quadratic growth in $L^1$
for the control via second-order conditions.

In the present paper,
we focus on the situation in which a Tikhonov term is present
(within the smooth part $F$)
and enables us to analyze the problem in $L^2$.
In this setting, we mainly reproduce the results of
\cite{CasasHerzogWachsmuth2015:1}.
However,
our approach has three advantages.
First, we were able to prove quadratic growth in an $L^2(\Omega\times(0,T))$-ball
in the case of $j = j_3$,
whereas the growth was only known to hold in an $L^\infty(\Omega;L^2(0,T))$-ball before.
This is important for the numerical analysis of such problems.
We note that optimality in an $L^2(\Omega\times(0,T))$-ball
has been shown
in \cite[Theorem~4.2]{CasasMateosRoesch2017}
for a problem without control constraints
and in \cite[Theorem~4.16]{CasasKunisch2024}
for a problem with infinite time horizon.
Second, our sufficient condition in the case $\bar u = 0$ and $j = j_2$
seems to be weaker.
More comments concerning these two points can be found in \cref{sec:application}.
Finally,
since we prove second-order epi-differentiability of the functionals $j_i$
(under mild assumptions),
the sensitivity analysis from \cite{ChristofWachsmuth2017:3}
is applicable to the problem at hand
and can be used to prove (directional) differentiability
of the solution w.r.t.\ possible perturbations of the data.
We also mention that we identified two issues with the analysis in
\cite{CasasHerzogWachsmuth2015:1}, see
\cref{lem:weak_jop}
and \cref{ex:critical_cone_j2}.

In order to analyze problem \eqref{eq:prob},
we use the reformulation
\begin{equation}
	\label{eq:problem}
	\text{Minimize}
	\quad
	F(u) + G_i(u)
	\qquad\text{w.r.t.\ } u \in L^2( \Omega \times (0,T) ),
\end{equation}
where
\begin{equation}
	\label{eq:G}
	G_i(u) = \delta_{\Uad}(u) + \mu j_i(u)
\end{equation}
with the indicator function $\delta_{\Uad} \colon L^2(\Omega \times (0,T)) \to \set{0,\infty}$
of the feasible set.
Since the functionals $G_i$
are nonsmooth (and even discontinuous everywhere),
it is not clear how
(directional) second-order derivatives should be defined.
In \cite{CasasHerzogWachsmuth2015:1},
the authors used an ad-hoc approach,
i.e., they defined reasonable expressions for the second-order derivatives
and proved that they can be used to arrive at second-order optimality conditions.
We follow the approach of
\cite{ChristofWachsmuth2017:1,WachsmuthWachsmuth2022}
and utilize the weak second-order subderivative of $G_i$
at $\bar u$ w.r.t.\ $w \in L^2(\Omega \times (0,T))$
defined for all directions $v \in L^2(\Omega)$
via
\begin{equation*}
	G_i''(\bar u, w; v)
	:=
	\inf
	\set*{
		\liminf_{k \to \infty} \frac{G_i(\bar u + t_k v_k) - G_i(\bar u) - t_k \dual{w}{v_k}}{t_k^2/2}
		\given
		t_k \searrow 0,
		v_k \weakly v
	}
	.
\end{equation*}
We prove that this weak second subderivative
coincides with the expressions given in
\cite{CasasHerzogWachsmuth2015:1}.
Further, we prove that the functionals $G_i$
are strongly twice epi-differentiable (see \cref{def:second_order_derivative})
and this enables us to
use abstract results concerning second-order optimality conditions.
Due to these preparations,
one can easily apply the theory to obtain
second-order conditions for problems in which, e.g.,
the functional $F$ is defined via different PDEs.

The paper is structured as follows.
In \cref{sec:second-order_conditions}
we give details concerning no-gap second-order conditions
via the calculus of subderivatives.
We review the second-order theory (\cref{subsec:known_second_order})
and also give some new results
of first and second order (\cref{subsec:new_first_order} and \cref{subsec:new_second_order}).
\Cref{sec:subderivatives_sparsity}
is devoted to the computation of
the weak second subderivatives of $G_i$
and to the verification of the strong twice epi-differentiability of $G_i$.
Finally, these findings are applied
to a semilinear parabolic control problem in \cref{sec:application}.

\section{No-gap second-order conditions}
\label{sec:second-order_conditions}
In this section, we consider the minimization problem
\begin{equation*}
	\label{eq:prob_abstract}
	\tag{P}
	\text{Minimize}\quad
	\Phi(x) := F(x) + G(x)
	\qquad\text{w.r.t. } x \in X.
\end{equation*}
Here, $G \colon X \to \barR := (-\infty, \infty]$
and $F \colon \dom(G) \to \R$
are given.
We are going to provide optimality conditions for
\eqref{eq:prob_abstract} by using subderivatives of $G$.

We are interested in necessary and sufficient conditions of second order,
such that the gap between both conditions is as small as in finite dimensions.
In \cref{subsec:known_second_order},
we present the second-order theory
from \cite{WachsmuthWachsmuth2022},
which is a slight generalization
of the theory from \cite{ChristofWachsmuth2017:1}.
Afterwards, we introduce a first-order subderivative
in \cref{subsec:new_first_order}
and provide associated results and calculus rules.
Finally,
in \cref{subsec:new_second_order},
we present some new results of second order.

Throughout this section, we always consider the following situation.
\begin{assumption}[Standing Assumptions and Notation]~
	\label{asm:standing_assumption}
	\begin{enumerate}
		\item
			\label{asm:standing_assumption:1}
			$X$ is the (topological) dual space of a separable Banach space $Y$,
		\item $\bar x \in \dom(G)$ is fixed,
		\item
			\label{asm:standing_assumption:3}
			There exist $F'(\bar x) \in Y$ and a bounded bilinear form $F''(\bar x) \colon X \times X \to \R$ with
			\begin{equation}
			\label{eq:hadamard_taylor_expansion}
				\lim_{k \to \infty}
				\frac{F(\bar x + t_k h_k) - F(\bar x) - t_k F'(\bar x) h_k - \frac12 t_k^2 F''(\bar x) h_k^2}{t_k^2}
				=
				0
			\end{equation}
			for all sequences $\seq{t_k} \subset \R^+ := (0,\infty)$, $\seq{h_k} \subset X$
			satisfying $t_k \searrow 0$,  $h_k \weaklystar h \in X$
			and
			$\bar x + t_k h_k \in \dom(G)$.
	\end{enumerate}
\end{assumption}

Note that we use the abbreviations $ F'(\bar x) h  := \dual{F'(\bar x)}{h}$
and $F''(\bar x) h^2 := F''(\bar x) ( h,h )$ for all $h \in X$ in \eqref{eq:hadamard_taylor_expansion},
and that \eqref{eq:hadamard_taylor_expansion} is automatically satisfied if $F$ admits a second-order Taylor expansion of the form
\begin{equation}
	\label{eq:strong_taylor_expansion}
	F(\bar x + h)
	-
	F(\bar x)
	-
	F'(\bar x) h
	-
	\frac12 F''(\bar x) h^2
	=
	\oo(\norm{h}^2_X)
	\quad\text{ as } \norm{h}_X \to 0.
\end{equation}

\subsection{Review of second-order theory}
\label{subsec:known_second_order}
First, we review the theory from
\cite{ChristofWachsmuth2017:1,WachsmuthWachsmuth2022}.
As a second derivative for the functional $G$,
we use the so-called
weak-$\star$ second subderivative.
\begin{definition}[Weak-\texorpdfstring{$\star$}{*} Second Subderivative]
	\label{def:weak_star_subderivative}
	Let $x \in \dom(G)$ and $w \in Y$ be given.
	The weak-$\star$ second subderivative
	$G''(x, w; \cdot) \colon X \to [-\infty, \infty]$
	of $G$ at $x$ for $w$ is defined via
	\begin{equation*}
		G''(x, w; h)
		:=
		\inf
		\set*{
			\liminf_{k \to \infty} \frac{G(x + t_k h_k) - G(x) - t_k \dual{w}{h_k}}{t_k^2/2}
			\given
			t_k \searrow 0,
			h_k \weaklystar h
		}
		.
	\end{equation*}
\end{definition}

In the case that $G$ is convex,
the parameter $w$
in $G''(x,w;\cdot)$
will often be taken from the convex subdifferential $\partial G(x) \subset X\dualspace$.
In the case that $X$ is not reflexive,
we identify $Y$ with a proper subspace of $X\dualspace = Y\bidualspace$ (in the canonical way).
In \cite[Remark~2.3]{WachsmuthWachsmuth2022}
an example is given
which shows that
the existence of $w\in Y \cap \partial G(x)$ is an additional regularity assumption,
since it requires the existence of subgradients of $G$ in the smaller pre-dual space $Y$.

We review some properties of $G''(x, w; \cdot )$.

\begin{lemma}[\texorpdfstring{\cite[Lemma~2.4]{WachsmuthWachsmuth2022}}{[Lemma~2.4, Wachsmuth and Wachsmuth, 2022]}]
	\label{lem:Gpp_convex}
	We assume that $G$ is convex and $x \in \dom(G)$.
	For $w \in Y \cap \partial G(x)$ we have
	\begin{align*}
		\forall h \in X : \qquad G''(x, w; h ) &\ge 0,
		\intertext{whereas in case $w \in Y \setminus \partial G(x)$ we have}
		\exists h \in X \setminus \set{0} : \qquad G''(x, w; h ) &= -\infty.
	\end{align*}
\end{lemma}

In the next definition, we ensure the existence of recovery sequences.
\begin{definition}[Second-Order Epi-Differentiability]
	\label{def:second_order_derivative}
	Let $x \in  \dom(G)$ and $w \in Y$ be given.
	The functional $G$ is said to be weak-$\star$ twice epi-differentiable
	(respectively, strictly twice epi-differentiable, respectively, strongly
	twice epi-differentiable)
	at $x$ for $w$ in a direction $h \in X$,
	if for all $\seq{t_k}\subset \R^+$ with $t_k \searrow 0$ there exists a sequence $\seq{h_k} \subset X$
	satisfying $h_k \weaklystar h$
	(respectively, $h_k \weaklystar h$ and $\norm{h_k}_X \to \norm{h}_X$, respectively, $h_k \to h$) and
	\begin{equation}
		\label{eq:recovery_sequence_def}
		G''(x, w; h )
		=
		\lim_{k \to \infty} \frac{G(x + t_k h_k) - G(x) - t_k \dual{w}{h_k}}{t_k^2/2}
		.
	\end{equation}
	The functional $G$ is called weak-$\star$/strictly/strongly twice epi-differentiable at $x$ for $w$ if it is
	weak-$\star$/strictly/strongly twice epi-differentiable at $x$ for $w$ in all directions $h \in X$.
\end{definition}
We note that a slightly weaker property
would be sufficient to apply the second-order theory below.
In fact, for every $h \in X$ we only need one pair of sequences $\seq{t_k}$ and $\seq{h_k}$
with the properties as in \cref{def:second_order_derivative},
see
\cite[Definition~3.4, Theorem~4.3]{ChristofWachsmuth2017:1}.
However, many functionals $G$ are actually
weak-$\star$/strictly/strongly
twice epi-differentiable
and this stronger property is also useful
for a differential sensitivity analysis,
see \cite{ChristofWachsmuth2017:3}.

Next, we state the second-order optimality conditions for \eqref{eq:prob_abstract}.
We start with the necessary condition
from \cite[Theorem~2.8]{WachsmuthWachsmuth2022}.

\begin{theorem}[Second-Order Necessary Condition]
	\label{thm:SNC}
	Suppose that $\bar x$ is a local minimizer of \eqref{eq:prob_abstract} such that
	\begin{equation}
		\label{eq:second_order_growth}
		\Phi(x)
		\ge
		\Phi(\bar x) + \frac{c}{2} \norm{x - \bar x}^2_X
		\qquad\forall x \in X, \norm{x - \bar x}_X \le \varepsilon
	\end{equation}
	holds for some $c \geq 0$ and some $\varepsilon > 0$.
	Assume further that one of the following conditions is satisfied.
	\begin{enumerate}
		\item The map $h \mapsto F''(\bar x) h^2$ is sequentially weak-$\star$ upper semicontinuous. \label{assumption-usc}
		\item The functional $G$ is strongly twice epi-differentiable at $\bar x$ for $-F'(\bar x)$. \label{assumption-mrc}
	\end{enumerate}
	Then
	\begin{equation}
		\label{eq:weakstarSNC}
		F''(\bar x) h^2
		+
		G''(\bar x, -F'(\bar x); h)
		\ge
		c \norm{h}_X^2
		\qquad
		\forall h \in X.
	\end{equation}
\end{theorem}

We continue with the sufficient condition,
see \cite[Theorem~2.9]{WachsmuthWachsmuth2022}.

\begin{theorem}[Second-Order Sufficient Condition]
	\label{thm:SSC_wo}
	Assume that the map $h \mapsto F''(\bar x) h^2$
	is sequentially weak-$\star$ lower semicontinuous
	and that
	\begin{equation}
		\label{eq:SSC_wo}
		F''(\bar x) h^2
		+
		G''(\bar x, -F'(\bar x); h)
		>
		0
		\qquad
		\forall h \in X \setminus \set{0}.
	\end{equation}
	Suppose further that
	\begin{equation*}
		\label{eq:NDC}
		\tag{\textup{NDC}}
		\begin{aligned}
			&\text{for all $\seq{t_k} \subset \R^+$, $\seq{h_k} \subset X$
			with $t_k \searrow 0$, $h_k \weaklystar 0$
			and $\norm{h_k}_X = 1$, we have}\\
			&\qquad
			\liminf_{k \to \infty} \parens[\bigg]{
				\frac1{t_k^2} \parens[\big]{G(\bar x + t_k h_k) - G(\bar x)}
				+
				\dual{F'(\bar x)}{h_k / t_k}
				+
				\frac12 F''(\bar x) h_k^2
			}
			> 0
			.
		\end{aligned}
	\end{equation*}
	Then $\bar x$ satisfies the growth condition \eqref{eq:second_order_growth} with some constants $c > 0$ and $\varepsilon > 0$.
\end{theorem}
In case $G$ is convex,
a first-order condition is actually hidden in \eqref{eq:SSC_wo},
see \cref{lem:Gpp_convex} and also \cref{cor:2nd_order_traditional} below.
The acronym \eqref{eq:NDC}
stands for ``non-degeneracy condition'', see \cite[Theorem~4.4]{ChristofWachsmuth2017:1}.
Sufficient conditions for \eqref{eq:NDC}
are given in \cite[Lemma~5.1]{ChristofWachsmuth2017:1}.
A slight generalization of their case (ii) applies to our
problem \eqref{eq:problem}.
\begin{lemma}
	\label{lem:Legendre_gives_NDC}
	Suppose that $G$ is convex, $X$ is a Hilbert space and $-F'(\bar x) \in \partial G(\bar x)$.
	We further require that
	$h \mapsto F''(\bar x) h^2$ is a Legendre form,
	i.e.,
	it is sequentially weakly lower semicontinuous and
	\begin{equation*}
		h_k \weakly h
		\quad\text{and}\quad
		F''(\bar x) h_k^2 \to F''(\bar x) h^2
		\qquad\Longrightarrow\qquad
		h_k \to h
	\end{equation*}
	holds for all sequences $\seq{h_k} \subset X$.
	Then, \eqref{eq:NDC} is satisfied.
\end{lemma}
\begin{proof}
	Let sequences $\seq{t_k}$ and $\seq{h_k}$ as in \eqref{eq:NDC} be given.
	Due to convexity, we have
	\begin{equation*}
		\frac1{t_k^2} \parens[\big]{G(\bar x + t_k h_k) - G(\bar x)}
		+
		\dual{F'(\bar x)}{h_k / t_k}
		\ge
		0.
	\end{equation*}
	Thus, it is sufficient to verify
	$ \liminf_{k \to \infty} F''(\bar x) h_k^2 > 0$.
	From the sequential weak lower semicontinuity
	and $h_k \weakly 0$,
	we already get
	$ \liminf_{k \to \infty} F''(\bar x) h_k^2 \ge 0$.
	The case
	$ \liminf_{k \to \infty} F''(\bar x) h_k^2 = 0$
	cannot appear,
	since this would lead to $h_k \to 0$
	(at least on a subsequence)
	in contradiction to $\norm{h_k}_X = 1$.
	This finishes the proof.
\end{proof}
We note that
\cite[Lemma~5.1(ii)]{ChristofWachsmuth2017:1}
is also formulated
for so-called Legendre-$\star$
forms on dual spaces of reflexive spaces
(see \cite[Definition~4.1.2]{Harder2021} for the terminology).
However, it was shown in
\cite[Theorem~4.3.9]{Harder2021}
that this setting already implies
that the underlying space is (isomorphic to)
a Hilbert space.

Finally,
we recall a result
which is helpful for the calculation
of second subderivatives
via dense subsets.
\begin{lemma}[\texorpdfstring{\cite[Lemma~3.2]{ChristofWachsmuth2017:3}}{[Lemma 3.2, Christof and Wachsmuth, 2020]}]
	\label{lem:density_lemma}
	Let $x \in \dom(G)$ and $w \in Y$ be given.
	We suppose the existence of a set $V \subset X$
	and a functional $Q\colon X \to [-\infty,\infty]$ such that
	\begin{enumerate}
		\item\label{lem:density_lemma:1}
		for all $h \in X$ we have $G''(x,w;h) \ge Q(h)$,
		\item\label{lem:density_lemma:2}
		for all $h \in V$ and all $\seq{t_k} \subset \R^+$ with $t_k \searrow 0$,
		there exists a sequence $\seq{h_k} \subset X$ satisfying
		$h_k \weaklystar h$, $\norm{h_k}_X \to \norm{h}_X$, and
		\begin{align*}
			Q(h) = \lim_{k \to \infty} \frac{G(x + t_k h_k) - G(x) - \dual{w}{h_k}}{t_k^2 / 2} \in [-\infty,\infty],
		\end{align*}
		\item\label{lem:density_lemma:3}
		for all $h \in X$
		with $Q(h) < \infty$
		there exists a sequence $\seq{h^l} \subset V$
		with $h^l \weaklystar h$, $\norm{h^l}_X \to \norm{h}_X$
		and $Q(h) \ge \liminf_{l \to \infty} Q(h^l)$.
	\end{enumerate}
	Then, $Q = G''(x,w;\cdot)$
	and $G$ is strictly twice epi-differentiable at $x$ for $w$.
	If, moreover, the sequences in
	\ref{lem:density_lemma:2} and \ref{lem:density_lemma:3}
	can be chosen strongly convergent,
	then $G$ is even strongly twice epi-differentiable at $x$
	for $w$.
\end{lemma}
Since we use a setting which is slightly different from
\cite[Lemma~3.2]{ChristofWachsmuth2017:3},
we give a sketch of the proof.
\begin{proof}
	Let $h \in X$ be given.
	From \ref{lem:density_lemma:1} and \ref{lem:density_lemma:2},
	we immediately get that $Q(h) = G''(x, w; h)$
	for all $h \in V$.

	From \ref{lem:density_lemma:1} we further get
	$G''(x, w; h) = \infty$
	if $Q(h) = \infty$.
	For such sequences, we can take $h_k \equiv h$ to obtain \eqref{eq:recovery_sequence_def}.
	Now, let $h \in X$ with $Q(h) < \infty$ and $\seq{t_k} \subset \R^+$ with $t_k \searrow 0$
	be arbitrary.
	Let $\seq{h^l} \subset V$ be the sequence from \ref{lem:density_lemma:3}.
	According to \ref{lem:density_lemma:2},
	we find sequences $\seq{h^l_k}_k \subset X$ with
	$h^l_k \weaklystar h^l$, $\norm{h^l_k}_X \to \norm{h^l}_X$,
	and
	\begin{align*}
		Q(h^l) = \lim_{k \to \infty} \frac{G(x + t_k h^l_k) - G(x) - \dual{w}{h^l_k}}{t_k^2 / 2} \in [-\infty,\infty].
	\end{align*}
	In order to handle the limits $\pm\infty$,
	we equip $[-\infty,\infty]$
	with the metric
	$\bar d(x_1, x_2) = \abs{\operatorname{atan}(x_1) - \operatorname{atan}(x_2)}$.
	Now we can continue as in
	the proof of \cite[Lemma~3.2]{ChristofWachsmuth2017:3},
	and create a diagonal sequence
	$\seq{\hat h_k} \subset V$
	such that $\hat h_k \weaklystar h$, $\norm{\hat h_k}_X \to \norm{h}_X$
	and
	\begin{equation*}
		Q(h)
		\ge
		\liminf_{l \to \infty} Q(h^l)
		\ge
		\lim_{k \to \infty} \frac{G(x + t_k \hat h_k) - G(x) - \dual{w}{\hat h_k}}{t_k^2 / 2}.
	\end{equation*}
	Finally,
	\ref{lem:density_lemma:1}
	allows to bound the left-hand side by
	$G''(x,w;h)$
	from above,
	while the right-hand side is bounded from below
	by this term via its definition.
	This shows that $\seq{\hat h_k}$ is a strictly convergent recovery sequence.
	Since $h \in X$ and $\seq{t_k} \subset \R^+$ were arbitrary,
	this shows that
	$G$ is strictly twice epi-differentiable at $x$ for $w$.

	The strong twice epi-differentiability can be proven along the same lines.
\end{proof}

\subsection{New results of first order}
\label{subsec:new_first_order}
In this section,
we define a subderivative of first order
and investigate its relations
with the second subderivative from \cref{def:weak_star_subderivative}.
We start with the definition.
\begin{definition}[Weak-\texorpdfstring{$\star$}{*} (First) Subderivative]
	\label{def:firstsubderivative}
	Let $x \in \dom(G)$ be given. The weak-$\star$ (first) subderivative of $G$ at $x$ in a direction $h \in X$ is defined by
	\begin{equation}
		\label{eq:lower_directional_epiderivative_def}
		G \firstsubderivative (x; h)
		:=
		\inf
		\set*{
			\liminf_{k \to \infty} \frac{G(x + t_k h_k) - G(x)}{t_k}
			\given
			t_k \searrow 0,
			h_k \weaklystar h
		}
		.
	\end{equation}
\end{definition}

We note that a similar derivative is the ``lower directional epiderivative''
in \cite[(2.68)]{BonnansShapiro2000}
if we apply this definition to the space $X$ equipped with its weak-$\star$ topology.
However, this definition cannot be stated by using sequences,
but one has to use weak-$\star$ convergent nets,
which is inconvenient, since weak-$\star$ convergent nets can be unbounded.
We further mention that the finite-dimensional analogue is a
classical object in variational analysis, see, e.g., \cite[Definition~8.1]{RockafellarWets1998}.

In \cite[Lemma~2.5]{WachsmuthWachsmuth2022}
it was shown that
the implication
\begin{equation*}
	G'(x; h) > \dual{w}{h}
	\;\Rightarrow\;
	G''(x, w; h) = + \infty
	\qquad\forall h \in X
\end{equation*}
holds
under the assumptions that $G$ is convex and that
the directional derivative $G'(x; \cdot)$
is sequentially weak-$\star$ lower semicontinuous.

By utilizing the weak-$\star$ first subderivative,
we are able to circumvent these additional assumptions.
The finite-dimensional analogue
was considered in \cite[(2.5)]{BenkoMehlitz2022}.
\begin{lemma}
	\label{lem:Gpp_directional_derivative}
	For all $x \in \dom(G)$ and $w \in Y$ we have
	\begin{align*}
		G \firstsubderivative (x; h) > \dual{w}{h}
		&\;\Rightarrow\;
		G''(x, w; h) = + \infty
		\qquad\forall h \in X,
		\\
		G \firstsubderivative (x; h) < \dual{w}{h}
		&\;\Rightarrow\;
		G''(x, w; h) = - \infty
		\qquad\forall h \in X.
	\end{align*}
\end{lemma}
\begin{proof}
	Let $h \in X$ with $G \firstsubderivative (x;h) > \dual{w}{h}$ be arbitrary.
	For all sequences
	$\seq{t_k} \subset \R^+$, $\seq{h_k} \subset X$
	with
	$h_k \weaklystar h$ and $t_k \searrow 0$ we have
	\begin{align*}
		\liminf_{k \to \infty} \frac{G(x + t_k h_k) - G(x) - t_k \dual{w}{h_k}}{t_k}
		&=
		\liminf_{k \to \infty} \frac{G(x + t_k h_k) - G(x)}{t_k}
		-
		\lim_{k \to \infty} \dual{w}{h_k}
		\\
		&\ge
		G \firstsubderivative (x;h) - \dual{w}{h}
		>
		0
		,
	\end{align*}
	where we used \cref{def:firstsubderivative}
	for ``$\ge$''.
	Since an additional factor $t_k^{-1}$ appears in the definition of $G''(x,w;h)$, this implies $G''(x,w;h) = + \infty$.
	The other implication follows by a similar argument.
\end{proof}

The next lemma
provides some properties of the subderivative.
In particular, it
shows that the subderivative coincides with the directional derivative
under certain assumptions on $G$.
Interestingly, these are the conditions appearing in \cite[Lemma~2.5]{WachsmuthWachsmuth2022}
and, therefore, this result
follows from \cref{lem:Gpp_directional_derivative,lem:directional_derivative_equals_first_epiderivative}.

\begin{lemma}
	\label{lem:directional_derivative_equals_first_epiderivative}
	Let $x \in \dom(G)$ be given.
	\begin{enumerate}
		\item\label{lem:directional_derivative_equals_first_epiderivative:1}
			Let a direction $h \in X$ be given, for which the directional derivative
			$G'(x; h) \in [-\infty,\infty]$ exists.
			Then,
			\begin{equation*}
				G'(x; h) \ge G\firstsubderivative(x; h).
			\end{equation*}
		\item\label{lem:directional_derivative_equals_first_epiderivative:4}
			If $G$ is convex, then
			$G\firstsubderivative(x; \cdot) \colon X \to [-\infty, \infty]$ is convex.
		\item\label{lem:directional_derivative_equals_first_epiderivative:2}
			If $G$ is convex, then
			for all $w \in Y$ we have
			\begin{equation*}
				w \in \partial G(x)
				\qquad\Leftrightarrow\qquad
				G\firstsubderivative(x; h)
				\ge
				\dual{w}{h}
				\quad\forall h \in X.
			\end{equation*}
		\item\label{lem:directional_derivative_equals_first_epiderivative:3}
			Suppose that
			$G$ is convex
			and that
			$G'(x; \cdot) : X \to [-\infty, \infty]$ is sequentially weak-$\star$ lower semicontinuous.
			Then,
			\begin{equation*}
				G'(x;h) = G \firstsubderivative (x;h)
				\qquad
				\forall h \in X.
			\end{equation*}
	\end{enumerate}
\end{lemma}
\begin{proof}
	\ref{lem:directional_derivative_equals_first_epiderivative:1}:
	The inequality
	$G'(x;h) \ge G \firstsubderivative (x;h)$
	follows directly from the definitions.

	\ref{lem:directional_derivative_equals_first_epiderivative:4}:
	Let $h_1, h_2 \in X$ and $\lambda \in (0,1)$ be given.
	We have to show
	\begin{equation*}
		G\firstsubderivative(x;  \lambda h_1 + (1 - \lambda) h_2 )
		\le
		\lambda G\firstsubderivative(x; h_1)
		+
		(1-\lambda) G\firstsubderivative(x; h_2),
	\end{equation*}
	where we use the convention $\infty + (-\infty) := (-\infty) + \infty := \infty$.
	Let $\seq{\tilde t_{j,k}}_k \subset \R^+$ and $\seq{\tilde h_{j,k}}_k \subset X$ with
	$\tilde t_{j,k} \searrow 0$
	and
	$\tilde h_{j,k} \weaklystar h_j$
	be arbitrary,
	where $j \in \set{1,2}$.
	We select subsequences, denoted by $\seq{t_{j,k}}_k$ and $\seq{h_{j,k}}_k$,
	such that
	\begin{equation*}
		\liminf_{k \to \infty} \frac{G(x + \tilde t_{j,k} \tilde h_{j,k}) - G(x)}{\tilde t_{j,k}}
		=
		\lim_{k \to \infty} \frac{G(x + t_{j,k} h_{j,k}) - G(x)}{t_{j,k}}
		\in
		[-\infty, \infty]
	\end{equation*}
	for $j \in \set{1,2}$,
	i.e., these subsequences realize the limit inferior.
	We set $t_k := \parens*{\lambda / t_{1,k} + (1-\lambda) / t_{2,k}}^{-1} \searrow 0$
	and $g_k := \lambda h_{1,k} + (1-\lambda) h_{2,k} \weaklystar \lambda h_1 + (1-\lambda) h_2$.
	By convexity of $G$ we get
	\begin{align*}
		G(x + t_k g_k)
		&=
		G\parens*{
			x
			+
			\frac{\lambda t_k}{t_{1,k}} (t_{1,k} h_{1,k})
			+
			\frac{(1-\lambda) t_k}{t_{2,k}} (t_{2,k} h_{2,k})
		}
		\\&
		\le
		\frac{\lambda t_k}{t_{1,k}} G(x + t_{1,k} h_{1,k})
		+
		\frac{(1-\lambda) t_k}{t_{2,k}} G(x + t_{2,k} h_{2,k})
		.
	\end{align*}
	Together with the definition of the subderivative, we find
	\begin{align*}
		\MoveEqLeft
		G\firstsubderivative(x;  \lambda h_1 + (1 - \lambda) h_2 )
		\\
		&\le
		\liminf_{k \to \infty} \frac{ G(x + t_k g_k) - G(x)}{t_k}
		\\
		&\le
		\liminf_{k \to \infty} \parens*{
			\lambda \frac{G(x + t_{1,k} h_{1,k}) - G(x)}{t_{1,k}}
			+
			(1-\lambda) \frac{G(x + t_{2,k} h_{2,k}) - G(x)}{t_{2,k}}
		}
		\\
		&\le
		\lambda \lim_{k \to \infty}
		\frac{G(x + t_{1,k} h_{1,k}) - G(x)}{t_{1,k}}
		+
		(1-\lambda) \lim_{k \to \infty}
		\frac{G(x + t_{2,k} h_{2,k}) - G(x)}{t_{2,k}}
		\\
		&=
		\lambda \liminf_{k \to \infty}
		\frac{G(x + \tilde t_{1,k} \tilde h_{1,k}) - G(x)}{\tilde t_{1,k}}
		+
		(1-\lambda) \liminf_{k \to \infty}
		\frac{G(x + \tilde t_{2,k} \tilde h_{2,k}) - G(x)}{\tilde t_{2,k}}
		.
	\end{align*}
	In the last inequality, the convention $\infty + (- \infty) = \infty$ is crucial.
	Since this holds for all sequences as above,
	we can take the infimum over these sequences and
	this yields the desired convexity.

	\ref{lem:directional_derivative_equals_first_epiderivative:2}:
	``$\Rightarrow$'':
	Let $\seq{t_k} \subset \R^+$ and $\seq{h_k} \subset X$
	with $t_k \searrow 0$ and $h_k \weaklystar h \in X$ be arbitrary.
	The definition of the subdifferential yields
	\begin{equation*}
		\frac{G(x + t_k h_k) - G(h)}{t_k}
		\ge
		\dual{w}{h_k}
		.
	\end{equation*}
	Due to $w \in Y$, we can pass to the limit inferior $k \to \infty$.
	Afterwards, we take the infimum over all such sequences
	$\seq{t_k}$, $\seq{h_k}$
	and this yields the claim.

	``$\Leftarrow$'':
	If
	$w \le G\firstsubderivative(x; \cdot)$,
	we get
	$w \le G'(x; \cdot)$
	from
	\ref{lem:directional_derivative_equals_first_epiderivative:1}
	and
	(together with the convexity of $G$)
	this implies
	$\dual{w}{h} \le G'(x; h) \le G(x + h) - G(x)$
	for arbitrary $h \in X$,
	i.e., $w \in \partial G(x)$.

	\ref{lem:directional_derivative_equals_first_epiderivative:3}:
	In view of
	\ref{lem:directional_derivative_equals_first_epiderivative:1},
	it is sufficient to prove
	``$\le$''.
	Let $\seq{t_k} \subset \R^+$ and $\seq{h_k} \subset X$
	with $t_k \searrow 0$ and $h_k \weaklystar h \in X$ be arbitrary.
	We get
	\begin{equation*}
		G'(x;h)
		\le
		\liminf_{k \to \infty} G'(x;h_k)
		\le
		\liminf_{k \to \infty} \frac{G(x + t_k h_k) - G(x)}{t_k}
		,
	\end{equation*}
	by the sequential weak-$\star$ lower semicontinuity of $G'(x; \cdot)$ and the convexity of $G$. Taking the infimum over all sequences yields the desired inequality.
\end{proof}

A simple argument leads to a first-order condition.
\begin{theorem}[First-Order Necessary Condition]
	\label{thm:FONC}
	Suppose that $\bar x$ is a local minimizer of \eqref{eq:prob_abstract}.
	Then,
	$F'(\bar x) h + G\firstsubderivative(\bar x; h) \ge 0$
	for all $h \in X$.
	If $G$ is additionally convex,
	this condition is equivalent to
	$-F'(\bar x) \in \partial G(\bar x)$.
\end{theorem}
\begin{proof}
	Let sequences $\seq{t_k} \subset \R^+$, $\seq{h_k} \subset X$
	with $t_k \searrow 0$ and $h_k \weaklystar h$ be given.
	From \eqref{eq:hadamard_taylor_expansion} we get
	\begin{equation*}
		\lim_{k \to \infty}
		\frac{F(\bar x + t_k h_k) - F(\bar x) - t_k F'(\bar x) h_k}{t_k}
		=
		0
		.
	\end{equation*}
	Thus,
	\begin{equation*}
		\liminf_{k \to \infty}
		\frac{G(\bar x + t_k h_k) - G(\bar x)}{t_k}
		\ge
		\liminf_{k \to \infty} \frac{F(\bar x) - F(\bar x + t_k h_k)}{t_k}
		=
		\lim_{k \to \infty} -F'(\bar x) h_k
		=
		-F'(\bar x) h.
	\end{equation*}
	The first part of the claim follows by taking the infimum
	w.r.t.\ all these sequences.
	The second claim follows from
	\itemref{lem:directional_derivative_equals_first_epiderivative:2}
	with $w = -F'(\bar x)$,
	since $F'(\bar x) \in Y$.
\end{proof}

Next,
we provide a sum rule for $G\firstsubderivative$.

\begin{lemma}
	\label{lem:sumrule_subderivative_take_2}
	Suppose that
	$G = g_1 + g_2$
	with $g_1, g_2 \colon X \to \barR$.
	For all $x \in \dom(G)$ we have
	\begin{equation*}
		G\firstsubderivative(x; h)
		\ge
		g_1\firstsubderivative(x; h)
		+
		g_2\firstsubderivative(x; h)
	\end{equation*}
	for all $h \in X$ for which the right-hand side
	is not $\infty + (-\infty)$ (or $(-\infty) + \infty$).
	Additionally, assume that $X$ is reflexive and
	\begin{enumerate}
		\item
			\label{lem:sumrule_subderivative_take_2: i}
			$g_1$ Lipschitz continuous in a neighborhood of $x$ and convex,
		\item
			\label{lem:sumrule_subderivative_take_2: ii}
			$g_2$ is strongly (once) epi-differentiable at $x$ in the sense
			that for all $h \in X$ and sequences $\seq{t_k} \subset \R^+$ with $t_k \searrow 0$
			there exists a sequence $\seq{h_k} \subset X$ with $h_k \to h$ in $X$
			such that
			\begin{equation*}
				g_2\firstsubderivative(x; h)
				=
				\lim_{k \to \infty}
				\frac{g_2(x + t_k h_k) - g_2(x)}{t_k}
				.
			\end{equation*}
	\end{enumerate}
	Then,
	\begin{equation*}
		G\firstsubderivative(x; h)
		=
		g_1\firstsubderivative(x; h)
		+
		g_2\firstsubderivative(x; h)
		\qquad\forall h \in X.
	\end{equation*}
\end{lemma}
\begin{proof}
	The inequality ``$\ge$'' follows by the definition of the first subderivative
	by using
	$\liminf_{n \to \infty} \parens{a_n + b_n} \ge \liminf_{n \to \infty} a_n + \liminf_{n \to \infty} b_n$,
	whenever
	the right-hand side is not $\infty + (-\infty)$.

	We show the other inequality under the additional assumptions.
	Let $\seq {t_k} \subset \R^+$ be arbitrary with $t_k \searrow 0$
	and let $\seq{h_k} \subset X$ with $h_k \to h$ be given according to
	\ref{lem:sumrule_subderivative_take_2: ii}.
	We have
	\begin{align*}
		\lim_{k \to \infty} \frac{g_1(x + t_k h_k) - g_1(x)}{t_k}
		&=
		\lim_{k \to \infty} \frac{g_1(x + t_k h) - g_1(x)}{t_k}
		+
		\lim_{\mathclap{k \to \infty}} \frac{g_1(x + t_k h_k) - g_1(x + t_k h)}{t_k}
		\\
		&=
		g_1'(x; h) + 0,
	\end{align*}
	where we used the convexity for the existence of the directional derivative
	and the Lipschitz continuity is applied for the second addend.
	Thus,
	\begin{align*}
		G \firstsubderivative (x; h)
		&\le
		\liminf_{k \to \infty} \frac{G(x + t_k h_k) - G(x)}{t_k} \\
		&=
		\lim_{k \to \infty} \frac{g_1(x + t_k h_k) - g_1(x)}{t_k} +
		\liminf_{k \to \infty} \frac{g_2(x + t_k h_k) - g_2(x)}{t_k}
		\\
		&=
		g_1'(x; h) + g_2 \firstsubderivative (x; h)
		.
	\end{align*}
	Finally, we note that
	$g_1'(x; \cdot)$
	is convex and Lipschitz continuous,
	thus sequentially weakly lower semicontinuous.
	Since $X$ is reflexive, this implies
	sequential weak-$\star$ lower semicontinuity
	and 
	\itemref{lem:directional_derivative_equals_first_epiderivative:3}
	implies
	$g_1'(x; h) = g_1\firstsubderivative(x; h)$
	for all $h \in X$.
\end{proof}

Finally,
we give two lemmas
in reflexive spaces.
The first of these results
is similar to \cite[Proposition~6.2]{MohammadiMordukhovich2021}.
\begin{lemma}
	\label{lem:subderivative_reflexive_1}
	We assume that the space $X$ is reflexive and
	that $G \colon X \to \barR$ is convex.
	For all $x \in \dom(G)$ and $h \in X$
	we have
	\begin{equation}
		\label{eq:firstsub_with_strong}
		G\firstsubderivative(x; h)
		=
		\inf
		\set*{
			\liminf_{k \to \infty} \frac{G(x + t_k h_k) - G(x)}{t_k}
			\given
			t_k \searrow 0,
			h_k \to h
		}
		.
	\end{equation}
	Moreover,
	there exist sequences 
	$\seq{t_k} \subset \R^+$ and $\seq{h_k} \subset X$ with $t_k \searrow 0$, $h_k \to h$
	and
	\begin{equation*}
		G\firstsubderivative(x; h)
		=
		\lim_{k \to \infty} \frac{G(x + t_k h_k) - G(x)}{t_k}
		.
	\end{equation*}
	In particular, $G\firstsubderivative(x; \cdot)$ is
	lower semicontinuous.
\end{lemma}
Note that \eqref{eq:firstsub_with_strong}
shows that the weak-$\star$ subderivative
coincides with the strong subderivative
(defined as in \cref{def:firstsubderivative} with strong convergence instead of weak-$\star$ convergence).
\begin{proof}
	We denote the right-hand side of \eqref{eq:firstsub_with_strong}
	by $R(h)$.
	It is clear that $R \ge G\firstsubderivative(x;\cdot)$.
	By arguing as in 
	\itemref{lem:directional_derivative_equals_first_epiderivative:4},
	we can check that $R$ is convex.
	Next, we verify that the infimum in the definition of $R$ is
	attained.
	For an arbitrary $h \in X$,
	the definition of $R(h)$ implies
	the existence of double sequences $\seq{t_{k,n}} \subset \R^+$, $\seq{h_{k,n}} \subset X$
	and
	$\seq{r_k} \subset [-\infty,\infty]$
	such that
	\begin{align*}
		\lim_{n \to \infty} \norm{h_{k,n} - h}_X &= 0 \qquad \forall k \in \N,
		&
		\lim_{n \to \infty} t_{k,n} &= 0 \qquad \forall k \in \N,
		\\
		\lim_{n \to \infty} \frac{G(x + t_{k,n} h_{k,n}) - G(x)}{t_{k,n}} &= \mrep{r_k}{0} \qquad \forall k \in \N,
		&
		\lim_{k \to \infty} r_k &= R(h).
	\end{align*}
	To handle the limits $\pm \infty$,
	we again use the metric
	$\bar d(x_1, x_2) = \abs{\operatorname{atan}(x_1) - \operatorname{atan}(x_2)}$
	on $[-\infty,\infty]$.
	For each $k \in \N$, we select $n(k) \in \N$ such that
	\begin{align*}
		\norm{h_{k,n(k)} - h}_X &\le \frac1k,
		&
		t_{k,n(k)} &\le \frac1k,
		&
		\bar d\parens[\bigg]{\frac{G(x + t_{k,n(k)} h_{k,n(k)}) - G(x)}{t_{k,n(k)}}, r_k} &\le \frac1k.
	\end{align*}
	This shows that the sequences
	$\seq{t_k} := \seq{t_{k,n(k)}}$,
	$\seq{h_k} := \seq{h_{k,n(k)}}$
	satisfy
	\begin{align*}
		h_k &\to h,
		&
		t_k &\searrow 0,
		&
		\frac{G(x + t_k h_k) - G(x)}{t_k} &\to R(h).
	\end{align*}
	Hence, the infimum in the definition of $R$
	is always attained.

	The lower semicontinuity of $R$
	follows by a similar diagonal-sequence argument.
	Together with the convexity,
	we get that $R$ is weakly lower semicontinuous.

	Now, let the sequences $\seq{t_k} \subset \R^+$, $\seq{h_k} \subset X$
	with $t_k \searrow 0$ and $h_k \weakly h$
	be arbitrary.
	From the convexity of $G$
	and the definition of $R$
	we get
	\begin{equation*}
		\frac{G(x + t_k h_k) - G(x)}{t_k}
		\ge
		\liminf_{s \searrow 0}\frac{G(x + s h_k) - G(x)}{s}
		\ge
		R(h_k).
	\end{equation*}
	The weak lower semicontinuity of $R$
	implies
	\begin{equation*}
		\liminf_{k \to \infty}
		\frac{G(x + t_k h_k) - G(x)}{t_k}
		\ge
		\liminf_{k \to \infty}
		R(h_k)
		\ge
		R(h).
	\end{equation*}
	Taking the infimum over all such sequences shows
	$G\firstsubderivative(x;h) \ge R(h)$ for all $h \in X$.
	This shows
	$G\firstsubderivative(x;\cdot) = R$
	and the claim follows.
\end{proof}
The reflexivity of $X$ is only used to
get the equivalence of
weak-$\star$ convergence
and weak convergence.
This is needed in order to apply the result
that convex and lower semicontinuous functionals
are weakly lower semicontinuous.
Without reflexivity,
we would get a similar assertion for the
weak subderivative of $G$
defined via
\begin{equation}
	\label{eq:weak_lower_directional_epiderivative_def}
	G \weakfirstsubderivative (x; h)
	=
	\inf
	\set*{
		\liminf_{k \to \infty} \frac{G(x + t_k h_k) - G(x)}{t_k}
		\given
		t_k \searrow 0,
		h_k \weakly h
	}
	.
\end{equation}
However,
if $X$ is not reflexive,
we cannot extract weakly convergent subsequences
and,
consequently,
$G \weakfirstsubderivative$
seems to be of little use.

\Cref{lem:subderivative_reflexive_1}
enables us to prove a very interesting
characterization of nonemptyness of the subdifferential.
\begin{lemma}
	\label{lem:subderivative_reflexive_2}
	Let $X$ be reflexive and $G \colon X \to \barR$ convex.
	For all $x \in \dom(G)$, the assertions
	\begin{enumerate}
		\item\label{lem:subderivative_reflexive_2:1}
			$G\firstsubderivative(x; 0) = 0$.
		\item\label{lem:subderivative_reflexive_2:3}
			$\partial G(x) \ne \emptyset$.
	\end{enumerate}
	are equivalent.
	In case that these assumptions hold, we also have
	\begin{equation}
		\label{eq:subderivative_via_subdifferential}
		G\firstsubderivative(x; h) = \sup\set{\dual{w}{h}_X \given w \in \partial G(x)}
		\qquad
		\forall h \in X
		.
	\end{equation}
\end{lemma}
\begin{proof}
	``\ref{lem:subderivative_reflexive_2:1}$\Rightarrow$\ref{lem:subderivative_reflexive_2:3}'':
	From \itemref{lem:directional_derivative_equals_first_epiderivative:4}
	and \cref{lem:subderivative_reflexive_1},
	we know that $G\firstsubderivative(x;\cdot)$
	is convex and lower semicontinuous.
	Further,
	for all $h \in X$
	we have
	\begin{equation*}
		0
		=
		G\firstsubderivative(x; 0)
		\le
		\liminf_{n \to \infty} G\firstsubderivative(x; \frac1n h)
		=
		\liminf_{n \to \infty} \frac1n G\firstsubderivative(x; h)
		=
		\begin{cases}
			+\infty & \text{if } G\firstsubderivative(x; h) = +\infty, \\
			0       & \text{if } G\firstsubderivative(x; h) \in \R   , \\
			-\infty & \text{if } G\firstsubderivative(x; h) = -\infty.
		\end{cases}
	\end{equation*}
	This
	implies that $G\firstsubderivative(x; \cdot)$
	never attains the value $-\infty$.
	Thus, we can invoke
	\cite[Proposition~I.3.1, p.~14]{EkelandTemam1999}
	and get
	\begin{equation*}
		G\firstsubderivative(x; h)
		=
		\sup\set[\big]{
			c + \dual{w}{h}_X
			\mid
			c \in \R, \; w \in Y, \;
			\forall \hat h \in X : c + \dual{w}{\hat h}_X \le G\firstsubderivative(x; \hat h )
		},
	\end{equation*}
	i.e., $G\firstsubderivative(x; \cdot)$
	is the pointwise supremum of its continuous, affine minorants.
	Since $G\firstsubderivative(x; \cdot)$ is positively homogeneous,
	one can check
	\begin{equation*}
		G\firstsubderivative(x; h)
		=
		\sup\set{
			\dual{w}{h}_X
			\mid
			w \in Y, \;
			\forall \hat h \in X : \dual{w}{\hat h}_X \le G\firstsubderivative(x; \hat h )
		}
		.
	\end{equation*}
	In order to conclude,
	it remains to check that
	$w \in \partial G(x)$
	if and only if
	$w \le G\firstsubderivative(x; \cdot)$,
	and this is precisely the assertion of
	\itemref{lem:directional_derivative_equals_first_epiderivative:2}.
	This shows
	\ref{lem:subderivative_reflexive_2:3}.
	Note that this part of the proof also shows validity of
	\eqref{eq:subderivative_via_subdifferential}.

	``\ref{lem:subderivative_reflexive_2:3}$\Rightarrow$\ref{lem:subderivative_reflexive_2:1}'':
	This follows directly from
	\itemref{lem:directional_derivative_equals_first_epiderivative:2}.
\end{proof}

One might wonder whether \eqref{eq:subderivative_via_subdifferential}
implies the conditions
\itemref{lem:subderivative_reflexive_2:1},
\ref{lem:subderivative_reflexive_2:3}.
It is clear that a linear and unbounded functional $G$
satisfies \eqref{eq:subderivative_via_subdifferential},
but
\itemref{lem:subderivative_reflexive_2:1},
\ref{lem:subderivative_reflexive_2:3}
are violated.
The next example shows that lower semicontinuity
also does not help.
It
provides a convex and lower semicontinuous function on a Hilbert space
such that its subderivative at $0$ is identically $-\infty$.
This example is similar to \cite[Remark~2.3]{WachsmuthWachsmuth2022}.

\begin{example}
	\label{ex:subderivative_not_proper}
	We consider the Hilbert space $\ell^2$ and the closed, convex set
	\begin{equation*}
		C := \set[\big]{ x \in \ell^2 \given \abs{x_n} \le n^{-2} \; \forall n \in \N}.
	\end{equation*}
	We define the function $f \colon \ell^2 \to \barR$ via
	\begin{equation*}
		f(x) :=
		\begin{cases}
			\sum_{n = 1}^\infty x_n & \text{for } x \in C \\
			+\infty & \text{else}.
		\end{cases}
	\end{equation*}
	It is clear that $f$ is convex.
	Since $n^{-2}$ is summable, we have $\dom(f) = C$
	and this set is closed.
	Next, we show that $f$ is continuous on its domain
	and this implies that $f$ is lower semicontinuous.
	To this end, let a sequence $\seq{x_m} \subset C$ be given
	such that $x_m \to x_0$ in $\ell^2$. Clearly, $x_0 \in C$.
	For an arbitrary $\varepsilon > 0$, there exists $N \in \N$
	with $\sum_{n = N+1}^\infty n^{-2} \le \varepsilon$.
	Due to $x_m \to x_0$, there exists $M \in \N$
	with $\sum_{n = 1}^N \abs{ (x_m - x_0)_n } \le \varepsilon$
	for all $m \ge M$.
	This implies
	\begin{equation*}
		\abs{ f(x_m) - f(x_0) }
		\le
		\sum_{n = 1}^N \abs{ (x_m - x_0)_n }
		+
		\sum_{n = N+1}^\infty \abs{ (x_m)_n } + \abs{ (x_0)_n }
		\le
		3 \varepsilon.
	\end{equation*}
	for all $m \ge M$.
	Since $\varepsilon > 0$ was arbitrary,
	this shows $f(x_m) \to f(x_0)$,
	i.e., $f$ is continuous on its domain.
	Together with the closedness of $\dom(f)$,
	the lower semicontinuity of $f$ follows.

	As announced, we demonstrate that
	$f\firstsubderivative(0; \cdot) \equiv -\infty$.
	We start by considering $h_0 \in c_c$,
	where $c_c \subset \ell^2$ is the (dense) subspace consisting of finite sequences.
	We choose $t_m := m^{-3} \searrow 0$
	and
	\begin{equation*}
		h_m := h_0 - \frac1m \sum_{n = m+1}^{m^2} e_n,
	\end{equation*}
	where $e_n$ is the $n$-th unit sequence.
	Note that $\norm{h_m - h_0}_{\ell^2}^2 = 1-1/m$
	and
	$\innerprod{h_m - h_0}{v}_{\ell^2} \to 0$
	for all $v \in c_c$.
	Thus,
	$h_m \weakly h_0$.
	For $m$ large enough, we have $t_m h_m \in C$.
	Thus,
	\begin{align*}
		f\firstsubderivative(0; h_0)
		\le
		\liminf_{m \to \infty} \frac{f(0 + t_m h_m) - f(0)}{t_m}
		&=
		\liminf_{m \to \infty} \sum_{n = 1}^\infty (h_m)_n
		\\&=
		\liminf_{m \to \infty} \sum_{n = 1}^\infty (h_0)_n - \frac{m^2 - m}{m}
		=
		-\infty
		.
	\end{align*}
	Thus, $f\firstsubderivative(0;\cdot)$
	equals $-\infty$ on the dense subspace $c_c$.
	Since $f\firstsubderivative(0;\cdot)$
	is lower semicontinuous due to \cref{lem:subderivative_reflexive_1},
	we get
	$f\firstsubderivative(0;\cdot) \equiv -\infty$.
	Finally, we can apply \cref{lem:subderivative_reflexive_2}
	to get $\partial f(0) = \emptyset$,
	although this can be also seen by elementary considerations.
	We also note that \eqref{eq:subderivative_via_subdifferential} still holds,
	since $\sup\emptyset = -\infty$.
\end{example}
A simple example shows that \eqref{eq:subderivative_via_subdifferential}
can also fail.
\begin{example}
	\label{ex:simple}
	Let $X = \R$ and consider $f \colon \R \to \barR$,
	\begin{equation*}
		f(x) :=
		\begin{cases}
			+\infty & \text{if } x < 0, \\
			-\sqrt{x} & \text{if } x \ge 0.
		\end{cases}
	\end{equation*}
	It is clear that $f$ is convex and lower semicontinuous.
	Further, we can check that
	\begin{equation*}
		f\firstsubderivative(0; h)
		=
		\begin{cases}
			+\infty & \text{if } h < 0, \\
			-\infty & \text{if } h \ge 0.
		\end{cases}
	\end{equation*}
	Thus, \itemref{lem:subderivative_reflexive_2:1}
	is violated.
	Consequently, \cref{lem:subderivative_reflexive_2}
	implies $\partial f(0) = \emptyset$.
	We also see that
	\eqref{eq:subderivative_via_subdifferential} fails.
\end{example}

\subsection{New results of second order}
\label{subsec:new_second_order}
Finally,
we present some new results
concerning
the second-order optimality conditions.

Surprisingly,
one can check that the quadratic growth condition
implies \eqref{eq:NDC}.
Note that this has not been realized in the previous works
\cite{ChristofWachsmuth2017:1,WachsmuthWachsmuth2022}.
\begin{theorem}[Quadratic Growth Implies \texorpdfstring{\eqref{eq:NDC}}{(NDC)}]
	\label{thm:NDC_necessary}
	We assume that the quadratic growth condition \eqref{eq:second_order_growth}
	is satisfied at $\bar x$ with some constants $c > 0$ and $\varepsilon > 0$.
	Then, \eqref{eq:NDC} is satisfied.
\end{theorem}
\begin{proof}
	Let sequences
	$\seq{t_k} \subset \R^+$, $\seq{h_k} \subset X$
	as in \eqref{eq:NDC} be given.
	From \eqref{eq:second_order_growth} we get
	\begin{align*}
		\frac c2
		&=
		\frac1{t_k^2} \parens*{ \frac c2 \norm{t_k h_k}_X^2 }
		\le
		\frac{F(\bar x + t_k h_k) - F(\bar x)}{t_k^2} + \frac{G(\bar x + t_k h_k) - G(\bar x)}{t_k^2}
	\end{align*}
	for $k$ large enough.
	Taking the limit inferior and using \eqref{eq:hadamard_taylor_expansion},
	we get
	\begin{align*}
		\frac c2
		&\le
		\liminf_{k \to \infty}
		\parens[\bigg]{
			\frac{F(\bar x + t_k h_k) - F(\bar x)}{t_k^2} + \frac{G(\bar x + t_k h_k) - G(\bar x)}{t_k^2}
		}
		\\
		&=
		\liminf_{k \to \infty}
		\parens[\bigg]{
			\frac{t_k F'(\bar x) h_k + \frac12 t_k^2 F''(\bar x) h_k^2 }{t_k^2}
			+ \frac{G(\bar x + t_k h_k) - G(\bar x)}{t_k^2}
		}
		\\
		&=
		\liminf_{k \to \infty} \parens[\bigg]{
			\frac1{t_k^2} \parens[\big]{G(\bar x + t_k h_k) - G(\bar x)}
			+
			\dual{F'(\bar x)}{h_k / t_k}
			+
			\frac12 F''(\bar x) h_k^2
		}.
	\end{align*}
	Since this holds for all sequences (as above),
	we obtain that \eqref{eq:NDC}
	is satisfied.
\end{proof}
This theorem even shows that the
limit inferior in \eqref{eq:NDC}
is uniformly positive.
It is also interesting to note that
this is always the case
whenever \eqref{eq:NDC}
is satisfied.
\begin{lemma}[Uniform Positivity in \texorpdfstring{\eqref{eq:NDC}}{(NDC)}]
	\label{lem:NDC_uniform}
	Suppose that
	\eqref{eq:NDC} is satisfied.
	Then, there exists $c > 0$
	such that
	\begin{equation*}
		\label{eq:NDC_uniform}
		\tag{\textup{NDC}'}
		\begin{aligned}
			&\text{for all $\seq{t_k} \subset \R^+$, $\seq{h_k} \subset X$
			with $t_k \searrow 0$, $h_k \weaklystar 0$
			and $\norm{h_k}_X = 1$, we have}\\
			&\qquad
			\liminf_{k \to \infty} \parens[\bigg]{
				\frac1{t_k^2} \parens[\big]{G(\bar x + t_k h_k) - G(\bar x)}
				+
				\dual{F'(\bar x)}{h_k / t_k}
				+
				\frac12 F''(\bar x) h_k^2
			}
			\ge
			\frac{c}{2}
		\end{aligned}
	\end{equation*}
	holds.
\end{lemma}
\begin{proof}
	We define
	\begin{equation*}
		\frac{c}{2}
		:=
		\inf\set*{
			\liminf_{k \to \infty} \parens[\bigg]{
				\frac{G(\bar x + t_k h_k) - G(\bar x)}{t_k^2}
				+
				\frac{\dual*{F'(\bar x)}{h_k}}{t_k}
				+
				\frac12 F''(\bar x) h_k^2
			}
			\given
			\begin{aligned}[c]
				&
				t_k \searrow 0,
				h_k \weaklystar 0,
				\\ &
				\norm{h_k}_X = 1
			\end{aligned}
		}.
	\end{equation*}
	In case $c = \infty$, there is nothing to show.
	Otherwise, we have $c \in [0,\infty)$.
	By definition,
	there are sequences of sequences
	$\seq{\seq{t_{k,n}}_n}_k$,
	$\seq{\seq{h_{k,n}}_n}_k$
	with $t_{k,n} \searrow 0$ and $h_{k,n} \weaklystar 0$
	as $n \to \infty$
	and $\norm{h_{k,n}}_X = 1$
	such that
	\begin{equation*}
		\frac{c_k}{2}
		:=
		\liminf_{n \to \infty} \parens[\bigg]{
			\frac{G(\bar x + t_{k,n} h_{k,n}) - G(\bar x)}{t_{k,n}^2}
			+
			\frac{\dual*{F'(\bar x)}{h_{k,n}}}{t_{k,n}}
			+
			\frac12 F''(\bar x) h_{k,n}^2
		}
	\end{equation*}
	satisfies $c_k \to c$.
	Now, we can argue as in the proof of
	\cite[Lemma~2.12(ii)]{ChristofWachsmuth2017:3}
	to obtain diagonal sequences
	$\seq{t_{k,n(k)}}$ and $\seq{h_{k,n(k)}}$
	with
	$t_{k,n(k)} \searrow 0$
	and
	$h_{k,n(k)} \weaklystar 0$
	as $k \to \infty$,
	$\norm{h_{k,n(k)}}_X = 1$
	and
	\begin{equation*}
		\lim_{k \to \infty} \parens[\bigg]{
			\frac{G(\bar x + t_{k,n(k)} h_{k,n(k)}) - G(\bar x)}{t_{k,n(k)}^2}
			+
			\frac{\dual[\big]{F'(\bar x)}{h_{k,n(k)}}}{t_{k,n(k)}}
			+
			\frac12 F''(\bar x) h_{k,n(k)}^2
		}
		=
		\frac{c}{2}
		.
	\end{equation*}
	From \eqref{eq:NDC},
	we infer $c > 0$
	and this yields the claim.
\end{proof}

By combining the previous theorem
with \cref{thm:SNC,thm:SSC_wo},
we arrive at our main theorem on no-gap second-order conditions.
\begin{theorem}[No-Gap Second-Order Optimality Condition]
	\label{thm:no-gap-SOC}
	Assume that the map $h \mapsto F''(\bar x) h^2$
	is sequentially weak-$\star$ lower semicontinuous
	and that one of the conditions \ref{assumption-usc} and \ref{assumption-mrc} in \cref{thm:SNC}
	is satisfied.
	Then,
	the second-order condition
	\eqref{eq:SSC_wo}
	and \eqref{eq:NDC} hold
	if and only if
	the
	quadratic growth condition \eqref{eq:second_order_growth} is satisfied
	at $\bar x$
	with constants $c>0$ and $\varepsilon > 0$.
\end{theorem}
We can recast the above second-order conditions
in a familiar form including
the first-order condition and a critical cone.
\begin{corollary}
	\label{cor:2nd_order_traditional_noncvx}
	Let the assumptions of \cref{thm:no-gap-SOC}
	and \eqref{eq:NDC} be satisfied.
	The quadratic growth condition \eqref{eq:second_order_growth}
	with constants $c > 0$ and $\varepsilon > 0$ is satisfied
	if and only if
	\begin{subequations}
		\label{eq:2nd_order_traditional_noncvx}
		\begin{align}
			\label{eq:2nd_order_traditional_first_order_noncvx}
			F'(\bar x) h + G\firstsubderivative(\bar x; h) &\ge 0 \qquad \forall h \in X, \\
			\label{eq:2nd_order_traditional_second_order_noncvx}
			F''(\bar x) h^2 + G''(\bar x, -F'(\bar x); h) &> 0 \qquad \forall h \in \KK \setminus \set{0},
		\end{align}
	\end{subequations}
	where
	the critical cone $\KK$ is defined via
	\begin{equation*}
		\KK := \set{ h \in X \given F'(\bar x) h + G \firstsubderivative (\bar x; h) = 0 }.
	\end{equation*}
	Moreover, if $\bar x$ is a local minimizer of \eqref{eq:prob_abstract},
	then \eqref{eq:2nd_order_traditional_noncvx} holds with
	``$\;\!\ge\!$'' instead of ``$\;\!>\!$'' in \eqref{eq:2nd_order_traditional_second_order_noncvx}.
\end{corollary}
\begin{proof}
	Let \eqref{eq:second_order_growth} be satisfied
	with $c \ge 0$ and $\varepsilon > 0$.
	From \cref{thm:SNC}, we get \eqref{eq:weakstarSNC}.
	Clearly, \eqref{eq:2nd_order_traditional_second_order_noncvx}
	(with ``$\ge$'' instead of ``$>$'' in case $c = 0$)
	follows.
	Since \eqref{eq:weakstarSNC}
	implies
	$G''(\bar x, -F'(\bar x); h) > -\infty$
	for all $h \in X$,
	\cref{lem:Gpp_directional_derivative}
	can be applied to obtain
	\eqref{eq:2nd_order_traditional_first_order_noncvx}.
	This shows the ``only if'' part of the first assertion
	and the second assertion.

	It remains to prove the ``if'' part of the first assertion.
	To this end,
	let \eqref{eq:2nd_order_traditional_noncvx} be satisfied.
	We show that this implies \eqref{eq:SSC_wo}.
	Let $h \in X \setminus \set{0}$ be given.
	In case $h \in \KK \setminus \set{0}$, \eqref{eq:SSC_wo}
	follows from \eqref{eq:2nd_order_traditional_second_order_noncvx}.
	Otherwise,
	$h \not\in \KK$
	and
	\eqref{eq:2nd_order_traditional_first_order_noncvx}
	give
	$F'(\bar x) h + G\firstsubderivative(\bar x; h ) > 0$.
	Thus,
	\eqref{eq:SSC_wo} is implied
	by \cref{lem:Gpp_directional_derivative}.
	Finally, \cref{thm:SSC_wo}
	shows that \eqref{eq:second_order_growth} holds.
\end{proof}

In the case that $G$ is convex,
we can use the assertions of \cref{lem:directional_derivative_equals_first_epiderivative}
to reformulate \eqref{eq:2nd_order_traditional_first_order_noncvx}
via the subdifferential.
\begin{corollary}
	\label{cor:2nd_order_traditional}
	In addition to the assumptions of \cref{thm:no-gap-SOC},
	we assume that \eqref{eq:NDC} holds and that $G$ is convex.
	The quadratic growth condition \eqref{eq:second_order_growth}
	with constants $c > 0$ and $\varepsilon > 0$ is satisfied
	if and only if
	\begin{subequations}
		\label{eq:2nd_order_traditional}
		\begin{align}
			\label{eq:2nd_order_traditional_first_order}
			F'(\bar x) + \partial G(\bar x) &\ni 0 \\
			\label{eq:2nd_order_traditional_second_order}
			F''(\bar x) h^2 + G''(\bar x, -F'(\bar x); h) &> 0 \qquad \forall h \in \KK \setminus \set{0},
		\end{align}
	\end{subequations}
	where
	the critical cone $\KK$ is defined via
	\begin{equation*}
		\KK := \set{ h \in X \given F'(\bar x) h + G \firstsubderivative (\bar x; h) = 0 }.
	\end{equation*}
	Moreover, if $\bar x$ is a local minimizer of \eqref{eq:prob_abstract},
	then \eqref{eq:2nd_order_traditional} holds with
	``$\;\!\ge\!$'' instead of ``$\;\!>\!$'' in \eqref{eq:2nd_order_traditional_second_order}.
\end{corollary}
\begin{proof}
	We just have to check that \eqref{eq:2nd_order_traditional_first_order_noncvx}
	and \eqref{eq:2nd_order_traditional_first_order} are equivalent,
	and this follows from
	\itemref{lem:directional_derivative_equals_first_epiderivative:2}.
\end{proof}

In the situation of \cref{cor:2nd_order_traditional},
let the first-order condition \eqref{eq:2nd_order_traditional_first_order}
be satisfied.
Using the equivalence with \eqref{eq:2nd_order_traditional_first_order_noncvx},
we get
$\KK = \set{ h \in X \given F'(\bar x) h + G\firstsubderivative(\bar x; h) \le 0 }$.
Since $G\firstsubderivative(\bar x; \cdot)$ is convex by
\itemref{lem:directional_derivative_equals_first_epiderivative:4},
this results in the convexity of $\KK$.
In case that $G\firstsubderivative(\bar x; \cdot)$ is additionally
(sequentially) weak-$\star$ lower semicontinuous,
$\KK$ is also (sequentially) weak-$\star$ closed.

For later reference,
we remark that the proofs of the last two corollaries
show that
\begin{equation}
	\label{eq:heile_welt_on_KK}
	G''(\bar x, -F'(\bar x); h) = +\infty
	\qquad\forall h \in X \setminus \KK
\end{equation}
holds whenever $G$ is convex and $-F'(\bar x) \in \partial G(\bar x)$.

Finally,
we also provide a sum rule for
$G''$.

\begin{lemma}[Sum Rule for Weak-\texorpdfstring{$\star$}{*} Second Subderivative]
	\label{lem:sumrule_second_subderivative}
	Let $g_1, g_2: X \to (-\infty, \infty]$ and $x \in \dom(g_1) \cap \dom(g_2)$.
	Furthermore $h \in X$ and $w_1, w_2 \in Y$.
	
	Then, it holds
	\begin{equation}
		(g_1 + g_2)''(x, w_1 + w_2; h) \ge g_1''(x, w_1; h) + g_2''(x, w_2; h),
	\end{equation}
	whenever the right-hand side is not $\infty + (-\infty)$.
\end{lemma}
\begin{proof}
	This follows from the definitions,
	see also the first part of the proof of \cref{lem:sumrule_subderivative_take_2}.
\end{proof}

\section{Second subderivatives of sparsity functionals}
\label{sec:subderivatives_sparsity}

Our plan is to apply the theory from \cref{sec:second-order_conditions}
to the problem \eqref{eq:problem}.
Throughout, we are using the spaces
$X = Y = \LtwoOmegaT$.
Here, $\OmegaT := \Omega \times (0,T)$,
where
$\Omega \subset \R^d$ is assumed to be non-empty, open, and bounded,
and $T > 0$.
Since the space $\LtwoOmegaT$ is reflexive,
the weak topology and the weak-$\star$ topology
coincide.
Thus, we will work with weak convergence throughout.
In particular, $v_k \weakly v$ always refers to weak convergence in $\LtwoOmegaT$
(unless stated otherwise).
Note that we will also write $\bar u \in \Uad$
instead of $\bar x$ to denote
the potential minimizer,
which is fixed throughout.

We use the functions $G_i$ from \eqref{eq:G},
i.e.,
\begin{equation*}
	G_i(u) = \mu j_i(u) + \delta_{\Uad}.
\end{equation*}
Here,
$\delta_{\Uad} \colon \LtwoOmegaT \to \set{0,\infty}$ is the indicator function (in the sense of convex analysis)
of the feasible set $\Uad$
and
$j_i$ is one of the functionals defined in
\eqref{eq:js},
scaled by $\mu > 0$.
Since the functions $j_i$ are finite everywhere on $\LtwoOmegaT$,
we get $\dom(G_i) = \Uad$.
At this point, we do not specify the function $F$,
we just require
that $F$ together with $F'(\bar u) \in \LtwoOmegaT$
and the bounded bilinear form $F''(\bar u) \colon \LtwoOmegaT \times \LtwoOmegaT \to \R$
satisfies
the Taylor-like expansion \eqref{eq:hadamard_taylor_expansion}.

We further recall that
the set of feasible controls is defined by
\begin{equation*}
	\Uad
	=
	\set{
		u \in \LtwoOmegaT
		\given
		\alpha \le u(x,t) \le \beta \text{ f.a.a.\ } (x,t) \in \Omega_T
	}
	,
\end{equation*}
where $\alpha,\beta \in \R$ are given with $\alpha < \beta$.
This set is convex, closed and bounded.
It is well known that the tangent cone (in the sense of convex analysis)
of $\Uad \subset \LtwoOmegaT$
at $\bar u \in \Uad$
is given by
\begin{equation}
	\label{eq:tangent_cone}
	\tangentcone{\bar u}
	= \set{
		v  \in \LtwoOmegaT
		\given
		v(x,t) \ge 0 \text{ if } \bar u(x,t) = \alpha
		\text{ and }
		v(x,t) \le 0 \text{ if } \bar u(x,t) = \beta
	}
\end{equation}
and concerning the normal cone
we have for all $v \in \LtwoOmegaT$ the equivalence
\begin{equation}
	\label{eq:normal_cone}
	v \in \normalcone{\bar u}
	\Leftrightarrow
	\begin{cases}
		v(x,t) \le 0 \text{ if } \bar u(x,t) = \alpha,
		\\
		v(x,t) \ge 0 \text{ if } \bar u(x,t) = \beta,
		\\
		v(x,t) = 0 \text{ if } \bar u(x,t) \in (\alpha, \beta)
	\end{cases}
	\text{f.a.a.\ $(x,t) \in \Omega_T$}.
\end{equation}
We expect that the results below can be extended to the case
in which the bounds $\alpha$ and $\beta$ are not constants.
For simplicity of the presentation,
we only consider the case of constant control bounds.

In all results of this section,
we do not use any special properties of
$(0,T)$ and $\Omega$.
Thus, they could be replaced by arbitrary
finite and complete measure spaces.
In particular, we could swap the roles of $(0,T)$
and $\Omega$.
This yields analogous results for the functionals
\begin{align*}
	j_4(u) &:= \norm{u}_{L^2(\Omega; L^1(0,T))} := \bracks*{ \int_\Omega \norm{u(x, \cdot)}_{L^1(0,T)}^2 \dx } ^ {1 / 2}, \\
	j_5(u) &:= \norm{u}_{L^1((0,T); L^2(\Omega))} := \int_0^T \norm{u(\cdot, t)}_{L^2(\Omega)} \dt
	,
\end{align*}
which involve sparsity w.r.t.\ time.

We start with the first-order analysis of problem \eqref{eq:problem}.
Since these preliminary results
hold for all $j_i$ or $G_i$, $i \in \set{1,2,3}$, we will just write $j$ or $G$.
\begin{lemma}
	It holds
	\begin{subequations}
		\begin{equation}
			\label{eq:first_subderivative_of_indicator_of_Uad}
			\delta_{\Uad}\firstsubderivative(\bar u; \cdot)
			=
			\delta_{\tangentcone{\bar u}}
		\end{equation}
		and
		\begin{equation}
			\label{eq:first_subderivative_G_form}
			G \firstsubderivative (\bar u; \cdot)
			=
			\mu j' (\bar u; \cdot) +
			\delta_{\tangentcone{\bar u}}
			.
		\end{equation}
	\end{subequations}
\end{lemma}
\begin{proof}
	In order to check \eqref{eq:first_subderivative_of_indicator_of_Uad},
	we use
	\cref{lem:subderivative_reflexive_1}
	and get
	\begin{equation*}
		\delta_{\Uad}\firstsubderivative(\bar u; v)
		=
		\inf
		\set*{
			\liminf_{k \to \infty} \frac{\delta_{\Uad}(\bar u + t_k v_k)}{t_k}
			\given
			t_k \searrow 0,
			v_k \to v
		}
		.
	\end{equation*}
	Now, it is straightforward
	to check that the right-hand side
	coincides with $\delta_{\tangentcone{\bar u}}(v)$.

	To show \eqref{eq:first_subderivative_G_form},
	we can apply \cref{lem:sumrule_subderivative_take_2} with
	$g_1 = \mu j$ and $g_2 = \delta_{\Uad}$.
	We verify the needed required properties.
	The space $\LtwoOmegaT$ is reflexive.
	\Itemref{lem:sumrule_subderivative_take_2: i} holds due to
	\begin{equation*}
		\abs{j_i(u_2) - j_i(u_1)}
		\le
		j_i(u_2 - u_1)
		\le
		C_i \Ltwonorm{u_2 - u_1}
		\qquad\forall u_1,u_2 \in \LtwoOmegaT.
	\end{equation*}
	For the convex and closed set $\Uad$,
	the Bouligand tangent cone $\tangentcone{\bar u}$
	coincides with the so-called inner tangent cone,
	see \cite[Proposition~2.55]{BonnansShapiro2000}.
	By the definition of the inner tangent cone,
	this result reads
	\begin{equation*}
		\tangentcone{\bar u}
		=
		\set*{
			v \in \LtwoOmegaT \given
			\forall \seq{t_k} \subset \R^+, t_k \searrow 0:
			\exists \seq{u_k} \subset \Uad :
			(u_k - \bar u) / t_k \to v
		}
	\end{equation*}
	and this is precisely property
	\itemref{lem:sumrule_subderivative_take_2: ii}.

	Applying \cref{lem:sumrule_subderivative_take_2} and
	\eqref{eq:first_subderivative_of_indicator_of_Uad} yields
	\begin{equation*}
		G \firstsubderivative (\bar u; v)
		=
		\mu j \firstsubderivative (\bar u; v) +
		\delta_{\Uad} \firstsubderivative (\bar u; v)
		=
		\mu j' (\bar u; v) +
		\delta_{\tangentcone{\bar u}}(v)
		,
	\end{equation*}
	where the last equality follows from
	\itemref{lem:directional_derivative_equals_first_epiderivative:3} as $j$ fulfills the required properties.
\end{proof}

Next, we see that the critical cone $\KK$ from \cref{cor:2nd_order_traditional}
coincides
with the critical cone $\criticalcone{\bar u}$
as defined in \cite[(4.1)]{CasasHerzogWachsmuth2015:1}.

\begin{lemma}
	\label{lem:critical_cones}
	Let $\bar u \in \Uad$ be given. Then
	the sets
	\begin{align*}
		\KK
		&:=
		\set{ v \in \LtwoOmegaT \given F'(\bar u) v + G \firstsubderivative (\bar u; v) = 0 },
		\\
		\criticalcone{\bar u}
		&:=
		\set {
			v \in \tangentcone{\bar u}
			\given
			F'(\bar u)v + \mu j'(\bar u; v) = 0
		}
	\end{align*}
	coincide.
\end{lemma}
\begin{proof}
	For $v \in \KK$,
	\eqref{eq:first_subderivative_G_form} yields
	$F'(\bar u)v + \mu j' (\bar u; v) + \delta_{\tangentcone{\bar u}}(v) = 0$.
	We have
	$\delta_{\tangentcone{\bar u}}(v) \in \set{0, \infty}$,
	but the value $\infty$
	would contradict the previous equality.
	This shows $v \in \criticalcone{\bar u}$.
	
	Now, let $v \in \criticalcone{\bar u}$ be given.
	From $v \in \tangentcone{\bar u}$
	and
	by using \eqref{eq:first_subderivative_G_form} again, we get
	\begin{equation*}
		0
		=
		F'(\bar u)v + \mu j'(\bar u; v)
		=
		F'(\bar u)v + \mu j'(\bar u; v) + \delta_{\tangentcone{\bar u}}(v)
		=
		F'(\bar u)v + G \firstsubderivative(\bar u; v)
		.
		\qedhere
	\end{equation*}
\end{proof}

We transfer the first-order necessary conditions
\eqref{eq:2nd_order_traditional_first_order}
for $\bar u \in \Uad$ to our situation.
	Since $j$ is continuous, the sum rule for the subdifferential applies,
	i.e.,
	\begin{equation*}
		\partial G(\bar u)
		=
		\mu \partial j(\bar u) + \partial \delta_{\Uad}(\bar u)
		=
		\mu \partial j(\bar u) + \normalcone{\bar u}.
	\end{equation*}
	Thus,
	the first-order condition
	\eqref{eq:2nd_order_traditional_first_order}
	is equivalent to
	the existence of
	$\lambda_{\bar u} \in \partial j(\bar u)$
	with
		\begin{equation}
			\label{necessary first order j 0}
			0 \in F'(\bar u) + \mu \lambda_{\bar u} + \normalcone{\bar u}
			.
		\end{equation}

Finally,
we characterize the directions from the critical cone.
\begin{lemma}
	\label{Lemma bei j1 ueber F + mu j' = 0}
	Suppose that
	$\lambda_{\bar u} \in \partial j(\bar u)$
	satisfies \eqref{necessary first order j 0}
	and let $v \in \tangentcone{\bar u}$
	be given.
	Then, $v \in \criticalcone{\bar u}$
	if and only if
	\begin{subequations}
		\label{Lemma j zu Zeigendes}
		\begin{align}
			\label{Lemma j zu Zeigendes:1}
			\parens{ F'(\bar u) + \mu \lambda_{\bar u} } v &= 0
			\qquad\text{a.e.\ in $\Omega_T$},
			\\
			j'(\bar u; v) &= \dual{\lambda_{\bar u}}{v}.
			\label{Lemma j zu Zeigendes:2}
		\end{align}
	\end{subequations}
\end{lemma}
\begin{proof}
	Let $v \in \criticalcone{\bar u}$ be given.
	With $w := \parens{ F'(\bar u) + \mu \lambda_{\bar u} } v \in L^1(\Omega_T)$
	we have
	\begin{equation*}
		0
		=
		F'(\bar u) v + \mu j'(\bar u; v)
		\ge
		\dual{F'(\bar u) + \mu \lambda_{\bar u} }{v}
		=
		\int_{\Omega_T} w \d(x,t)
		\ge
		0
	\end{equation*}
	by definition of $\criticalcone{\bar u}$,
	the properties of the subdifferential,
	and
	\eqref{necessary first order j 0}
	as
	$v \in \tangentcone{\bar u}$.
	This shows \eqref{Lemma j zu Zeigendes:2}
	and
	that the integral over $w$
	is zero.
	From the characterizations of the tangent cone and the normal cone,
	we get $w \ge 0$
	a.e.\ in $\Omega_T$.
	Thus, \eqref{Lemma j zu Zeigendes:1} follows.

	The converse implication follows immediately.
\end{proof}

\subsection{Second subderivative of \texorpdfstring{$j_1$}{j1}}

At the beginning, we recall the subdifferential of $j_1$.
For a functional $\lambda \in \LtwoOmegaT$
we have $\lambda \in \partial j_1(\bar u)$ if and only if
\begin{equation}\label{eq:subdiff_j1}
	\lambda(x,t) \in
	\Sign(\bar u(x,t))
	\qquad
	\text{f.a.a.\ } (x,t) \in \Omega_T
	,
\end{equation}
with the set-valued signum function
\begin{equation}
	\label{eq:Sign}
	\Sign (\theta) :=
	\begin{cases}
		\set{+1}					&\text{if } \theta > 0, \\
		\set{-1}					&\text{if } \theta < 0, \\
		[-1,+1]	&\text{if } \theta = 0.
	\end{cases}
\end{equation}
Moreover,
the directional derivative $j_1(\bar u;\cdot) : \LtwoOmegaT \to \R$ is given by
\begin{equation}
	\label{eq:directional_derivative_j1}
	j_1'(\bar u;v)
	=
	\int_{\set{\bar u > 0}} v \dxt
	- \int_{\set{\bar u < 0}}{v} \dxt
	+ \int_{\set{\bar u = 0}}{\abs{v}} \dxt.
\end{equation}
In order to apply the characterization of the critical cone
from \cref{Lemma bei j1 ueber F + mu j' = 0},
we analyze condition \eqref{Lemma j zu Zeigendes:2}
for $j = j_1$.

\begin{lemma}
	\label{lem:equality_subdif_directional_1}
	Let $\lambda_{\bar u} \in \partial j_1(\bar u)$
	be given.
	For $v \in \LtwoOmegaT$,
	we have
	$\dual{\lambda_{\bar u}}{v} = j_1'(\bar u; v)$
	if and only if
	\begin{equation*}
		\lambda_{\bar u}(x,t) =
		\begin{cases}
			+1, &\text{if } \bar u(x,t) = 0 \text{ and } v(x,t) > 0 \\
			-1, &\text{if } \bar u(x,t) = 0 \text{ and } v(x,t) < 0
		\end{cases}
		\qquad
		\text{f.a.a.\ } (x,t) \in \Omega_T.
	\end{equation*}
\end{lemma}
\begin{proof}
	This follows easily from
	\begin{equation*}
		j_1'(\bar u; v) - \dual{\lambda_{\bar u}}{v}
		=
		\int_{\set{\bar u = 0}} \abs{v} - \lambda_{\bar u} v \d(x,t),
	\end{equation*}
	since the integrand is non-negative a.e.\ due to \eqref{eq:subdiff_j1}.
\end{proof}
For the verification of the strong twice epidifferentiability,
we use the recovery sequence
from \cite[Theorem~3.7]{CasasHerzogWachsmuth2012}.
\begin{lemma}\label{lem:special_vk_for_j1}
	Assume $-F'(\bar u) \in \partial G_1(\bar u)$
	and let $v \in \criticalcone{\bar u}$ be given.
	Furthermore, let $\seq{t_k} \subset \R^+$ be an arbitrary sequence with $t_k \searrow 0$.
	We define
	the sequence $\seq{v_k} \subset \LtwoOmegaT$
	via
	(pointwise)
	\begin{equation}
		\label{eq:special_vk_for_j1}
		v_k :=
		\begin{cases}
			0                             &\text{if } \bar u \in (\alpha, \alpha +\sqrt{t_k}) \cup (\beta - \sqrt{t_k}, \beta) \cup (-\sqrt{t_k}, 0) \cup (0, \sqrt{t_k}), \\
			P_{t_k}(v) &\text{otherwise},
		\end{cases}
	\end{equation}
		where $P_{t_k} \colon \R \to \R$ denotes the projection onto the interval
			$\bracks*{ -\frac 1 {\sqrt{t_k}}, \frac 1 {\sqrt{t_k}} }$.
	Then, we have
	\begin{subequations}
		\label{eq:special_vk_for_j1_properties}
		\begin{align}
			\label{eq:special_vk_for_j1:1}
			v_k &\to v \quad\text{in } \LtwoOmegaT, \\
			\label{eq:special_vk_for_j1:2}
			v_k &\in \criticalcone{\bar u}, \\
			\label{eq:special_vk_for_j1:3}
			\bar u + t_k v_k &\in \Uad, \\
			\label{eq:special_vk_for_j1:4}
			j_1(\bar u + t_k v_k) - j_1(\bar u) &= t_k j_1'(\bar u; v_k).
		\end{align}
	\end{subequations}
\end{lemma}
\begin{proof}
	It holds $v_k(x,t) \to v(x,t)$ pointwise
	and $\abs{v_k(x,t)} \le \abs{v(x,t)}$.
	Due to $v \in \LtwoOmegaT$,
	Lebesgue's dominated convergence theorem yields \eqref{eq:special_vk_for_j1:1}.
	
	By construction, we get \eqref{eq:special_vk_for_j1:3}.
	
	To show \eqref{eq:special_vk_for_j1:2} we first note that
	$v_k \in \tangentcone{\bar u}$
	follows from \eqref{eq:special_vk_for_j1:3}.
	Next, we fix $\lambda_{\bar u} \in \partial j(\bar u)$
	such that \eqref{necessary first order j 0} holds.
	From \cref{Lemma bei j1 ueber F + mu j' = 0},
	we get that $v$ satisfies \eqref{Lemma j zu Zeigendes}.
	Using
	$\set{v_k \ne 0} \subset \set{v \ne 0}$
	by construction,
	we immediately get that \eqref{Lemma j zu Zeigendes:1}
	holds with $v$ replaced by $v_k$.
	Next, \eqref{Lemma j zu Zeigendes:2} enables us to apply \cref{lem:equality_subdif_directional_1}
	to $v$.
	Due to
	$\set{v_k > 0} \subset \set{v > 0}$
	and
	$\set{v_k < 0} \subset \set{v < 0}$,
	we can consequently invoke \cref{lem:equality_subdif_directional_1}
	with $v_k$ to obtain
	$\dual{\lambda_{\bar u}}{v_k} = j_1'(\bar u; v_k)$.
	Thus,
	\cref{Lemma bei j1 ueber F + mu j' = 0} can be applied to $v_k$
	to get $v_k \in \criticalcone{\bar u}$.

	Lastly, we prove \eqref{eq:special_vk_for_j1:4}.
	Easy calculations show that $\bar u(x,t) > 0$ implies
	$\bar u(x,t) + t_k v_k(x,t) \ge 0$.
	Analogously, $\bar u(x,t) + t_k v_k(x,t) < 0$ holds whenever $\bar u(x,t) < 0$.
	This yields
	\begin{align*}
		j_1(\bar u + t_k v_k) - j_1(\bar u)
		 &= \int_{\Omega_T} {\abs{\bar u + t_k v_k} - \abs{\bar u}} \dxt
		 \\
		 &=
		 t_k \parens*{
		 	\int_{\set{\bar u > 0}} {v_k} \dxt
		 	- \int_{\set{\bar u < 0}} {v_k} \dxt
			+ \int_{\set{\bar u = 0}} {\abs*{ v_k }} \dxt
		}
		 ,
	\end{align*}
	which shows the claim.
\end{proof}

\begin{theorem}\label{G1''}
	We assume
	$-F'(\bar u) \in \partial G_1(\bar u)$.
	Then,
	\begin{equation*}
		G_1''(\bar u,-F'(\bar u);v) = 0
		\qquad\forall v \in \criticalcone{\bar u}
	\end{equation*}
	holds and $G_1$ is strongly twice epi-differentiable at $\bar u$ for $-F'(\bar u)$.
\end{theorem}
\begin{proof}
	\Cref{lem:Gpp_convex} yields $G_1''(\bar u,-F'(\bar u);v) \ge 0$.
	Now, let $v \in \criticalcone{\bar u}$ and $\seq{t_k} \subset \R^+$ with $t_k \searrow 0$ be arbitrary.
	For the sequence $\seq{v_k} \subset \LtwoOmegaT$
	defined in \eqref{eq:special_vk_for_j1}
	we get
	\begin{align*}
		\MoveEqLeft
		\lim_{k \to \infty}
		\frac{G_1(\bar u + t_k v_k) - G_1(\bar u) - t_k \dual{-F'(\bar u)}{v_k}}{t_k ^2 / 2} \\
		&=
		\lim_{k \to \infty}
		\frac 2 {t_k^2} \parens*{ \indicatorofUad(\bar u + t_k v_k) + \mu j_1(\bar u + t_k v_k) - 0 - \mu j_1(\bar u) + t_k \dual{F'(\bar u)}{v_k} }
		\\
		&=
		\lim_{k \to \infty} \frac 2 {t_k} \parens*{ \mu j_1'(\bar u;v_k) + \dual{F'(\bar u)}{v_k} }
		&\qquad&\mathllap{\text{(by \eqref{eq:special_vk_for_j1:3} and \eqref{eq:special_vk_for_j1:4})}}\\
		&=
		\lim_{k \to \infty} \frac 2 {t_k} 0 = 0.
		&&\mathllap{\text{(by \eqref{eq:special_vk_for_j1:2})}}
	\end{align*}
	This shows $G_1''(\bar u,-F'(\bar u);v) = 0$
	and also the
	strong twice epi-differentiability
	in direction $v \in \criticalcone{\bar u}$.
	For $v \in \LtwoOmegaT \setminus \criticalcone{\bar u}$,
	we have $G_1''(\bar u, -F'(\bar u); v) = \infty$,
	see \eqref{eq:heile_welt_on_KK},
	and, hence, we can use $v_k \equiv v$
	as a recovery sequence.
\end{proof}

\subsection{Second subderivative of \texorpdfstring{$j_2$}{j2}}
As in
\cite[Section~3.2]{CasasHerzogWachsmuth2015:1},
we define $j_\Omega : L^2(\Omega) \to \R$ via
\begin{equation*}
	j_\Omega(u) := \normLoneOmega{u} = \int_{\Omega}{\abs*{ u(x) }} \dx.
\end{equation*}
The directional derivative $j_\Omega'(u; \cdot) : L^2(\Omega) \to \R$ is given by
\begin{equation}
	\label{jOmegap}
	j_\Omega'(u;v) = \int_{\set{u > 0}}{v(x)} \dx - \int_{\set{u < 0}}{v(x)} \dx + \int_{\set{u = 0}}{\abs{v(x)}} \dx.
\end{equation}
Now we can write $j_2$ as
\begin{equation}
	\label{eq: other form of j2}
	j_2(u) = \parens[\bigg]{ \int_0^T {j_\Omega(u(t))^2} \dt } ^{1/2}.
\end{equation}
The directional derivative $j_2'(u;\cdot) : \LtwoOmegaT \to \R$ is given by
\begin{equation}
	\label{eq:directional_derivative_j2}
	j_2'(u;v) =
	\begin{cases}
		j_2(v),
		&\text{if } u = 0,
		\\
		\frac 1 {j_2(u)} \int_0^T {j_\Omega'(u(t);v(t)) \normLoneOmega{u(t)}} \dt,
		&\text{if } u \ne 0,
	\end{cases}
\end{equation}
see \cite[Proposition~3.5]{CasasHerzogWachsmuth2015:1}.
Furthermore, it holds $\lambda \in \partial j_2(u)$ if and only if $\lambda \in L^2(0,T;L^\infty(\Omega))$ and
\begin{equation}
	\label{eq:subdiff_j2}
	\begin{split}
		u \ne 0: &\quad 
		\lambda(x,t) \in \Sign(u(x,t)) \frac {\normLoneOmega{u(t)}} {j_2(u)} \text{ a.e.\ in } \Omega_T
		\\
		u = 0: &\quad
		\norm{\lambda}_{L^2(0,T;L^\infty(\Omega))} \le 1, \\
	\end{split}
\end{equation}
see \cite[Proposition~3.4]{CasasHerzogWachsmuth2015:1}.
Note that $L^2(0,T;L^\infty(\Omega))$
is not a Bochner--Lebesgue space,
but the (canonical) dual of $L^2(0,T; L^1(\Omega))$,
which consists of weak-$\star$ measurable functions.
Here, we used the set-valued signum function
from \eqref{eq:Sign}.
Next, we recall a lower Taylor expansion of $j_2(\bar u + v)$.

\begin{lemma}
	[\texorpdfstring{\cite[Lemma~5.7]{CasasHerzogWachsmuth2015:1}}{[Lemma 5.7, Casas et al., 2017]}]
	\label{lem:lower_taylor_expansion}
	We assume $\bar u \ne 0$. There exist $C > 0$ and $\varepsilon > 0$ such that
	\begin{equation}
		\label{taylor estimate for j2} 
		j_2(\bar u + v)
		\ge
		j_2(\bar u)
		+
		j_2'(\bar u;v)
		+
		\frac 1 {2 j_2(\bar u)} \braces*{ \int_0^T {{j_\Omega'(\bar u(t);v_k(t))}^2} \dt - j_2'(\bar u;v_k)^2 }
		-
		C \frac {\Ltwonorm{v}^3} {j_2(\bar u)^2}
	\end{equation}
	holds for all $\Ltwonorm{v} \le \varepsilon$.
\end{lemma}

We provide a corrected version of
\cite[Lemma~5.6]{CasasHerzogWachsmuth2015:1}.
\begin{lemma}
	\label{lem:weak_jop}
	Let $\seq{v_k} \subset L^2(\Omega_T)$ be a sequence
	with $v_k \weakly v$ in $L^2(\Omega_T)$
	and $j_2'(\bar u; v_k) \to j_2'(\bar u;v)$.
	Then,
	the functions
	$w_k \in L^2(0,T)$,
	defined via
	$w_k(t) = \chi_M(t) j_\Omega'(\bar u(t), v_k(t))$,
	converge weakly in $L^2(0,T)$
	towards
	$w \in L^2(0,T)$,
	$w(t) = \chi_M(t) j_\Omega'(\bar u(t), v(t))$,
	where
	$M := \set{ t \in (0,T) \given \norm{\bar u(t)}_{L^1(\Omega)} \ne 0 }$.
\end{lemma}
\begin{proof}
	It is sufficient to consider the case $\bar u \ne 0$
	since $M = \emptyset$ if $\bar u = 0$.
	We define the functions
	\begin{align*}
		\hat w_k &:= \chi_{N_+} v_k - \chi_{N_-} v_k + \chi_{N_0} \abs{v_k}, &
		\hat w   &:= \chi_{N_+} v   - \chi_{N_-} v   + \chi_{N_0} \abs{v  },
	\end{align*}
	with the sets
	\begin{align*}
		N_+ &:= \set{(x,t) \in \Omega_T \given \norm{\bar u(t)}_{L^1(\Omega)} \ne 0, \; \bar u(x,t) > 0}, \\
		N_0 &:= \set{(x,t) \in \Omega_T \given \norm{\bar u(t)}_{L^1(\Omega)} \ne 0, \; \bar u(x,t) = 0}, \\
		N_- &:= \set{(x,t) \in \Omega_T \given \norm{\bar u(t)}_{L^1(\Omega)} \ne 0, \; \bar u(x,t) < 0}.
	\end{align*}
	Our first goal is to check $\hat w_k \weakly \hat w$ in $L^2(\Omega_T)$.
	Clearly, we already have weak convergence in the first two addends
	and it remains to consider $\chi_{N_0} \abs{v_k}$.
	Since this sequence is bounded in the reflexive space $L^2(\Omega_T)$,
	we get weak convergence of a subsequence (without relabeling),
	i.e.,
	$\chi_{N_0} \abs{v_k} \weakly z$ in $L^2(\Omega_T)$.
	For an arbitrary measurable set $Q \subset \Omega_T$, this yields
	\begin{equation*}
		\int_Q z \dxt
		=
		\lim_{k \to \infty} \int_Q \chi_{N_0} \abs{v_k} \dxt
		\ge
		\lim_{k \to \infty} \abs*{\int_Q \chi_{N_0} v_k \dxt}
		=
		\abs*{\int_Q \chi_{N_0} v \dxt}
		.
	\end{equation*}
	Consequently,
	$z \ge \chi_{N_0} \abs{v}$ a.e.\ on $\Omega_T$.
	Utilizing the formula for $j_2'$ and the assumption,
	we get
	\begin{align*}
		\lim_{k \to \infty} j_2'(\bar u; v_k)
		=
		j_2'(\bar u; v)
		&=
		\frac1{j_2(\bar u)}
		\int_{\Omega_T} \parens*{\chi_{N_+} v - \chi_{N_-} v + \chi_{N_0} \abs{v} } \norm{\bar u(t)}_{L^1(\Omega)} \dxt
		.
	\end{align*}
	On the other hand,
	\begin{align*}
		\lim_{k \to \infty} j_2'(\bar u; v_k)
		&=
		\lim_{k \to \infty}
		\frac1{j_2(\bar u)}
		\int_{\Omega_T} \parens*{\chi_{N_+} v_k - \chi_{N_-} v_k + \chi_{N_0} \abs{v_k}} \norm{\bar u(t)}_{L^1(\Omega)} \dxt
		\\
		&=
		\frac1{j_2(\bar u)}
		\int_{\Omega_T} \parens*{\chi_{N_+} v - \chi_{N_-} v + z} \norm{\bar u(t)}_{L^1(\Omega)} \dxt
		.
	\end{align*}
	Since the limit is unique, we get
	\begin{equation*}
		0
		=
		\int_{\Omega_T} \parens*{z - \chi_{N_0} \abs{v} } \norm{\bar u(t)}_{L^1(\Omega)} \dxt.
	\end{equation*}
	Since the integrand is nonnegative,
	it has to vanish a.e.\ on $\Omega_T$.
	The function $z$ satisfies $z = 0$ a.e.\ on $\Omega_T \setminus N_0$
	and we have $\norm{\bar u(t)}_{L^1(\Omega)} > 0$ for a.a.\ $(x,t) \in N_0$.
	Hence, $z = \chi_{N_0} \abs{v}$ a.e.\ on $\Omega_T$.
	This shows that the weak limit $z$ of $\chi_{N_0} \abs{v_k}$ is uniquely determined,
	consequently,
	the usual subsequence-subsequence argument yields the convergence of the entire sequence.
	Thus, have shown that
	$\hat w_k \weakly \hat w$ in $L^2(\Omega_T)$.
	
	Finally, the sequence $w_k$ is just the image of $\hat w_k$ under the bounded linear mapping
	$\int_\Omega \cdot \dx \colon L^2(\Omega_T) \to L^2(0,T)$, i.e.,
	\(
		w_k(t) = \int_\Omega \hat w_k(x,t) \dx.
	\)
	Thus, we get the desired $w_k \weakly w$ in $L^2(0,T)$.
\end{proof}
Let us briefly comment on the flaw in
\cite[Lemma~5.6]{CasasHerzogWachsmuth2015:1}.
Therein, the assertion of \cref{lem:weak_jop}
was proved for $M$ being the entire interval $(0,T)$.
This cannot be true
since the assumptions do not contain any information on
$v_k(\cdot, t)$ if $\norm{\bar u(t)}_{L^1(\Omega)} = 0$,
see \eqref{eq:directional_derivative_j2}.
Concerning the proof,
note that
\cite[(5.24)]{CasasHerzogWachsmuth2015:1}
reads ``$0 - 0 \to 0$'' for all $(x,t)$ with $\norm{\bar u(t)}_{L^1(\Omega)} = 0$,
but afterwards, the authors divide by $0$.
Finally, we mention that
\cite[Lemma~5.6]{CasasHerzogWachsmuth2015:1}
is only used in the proof of
\cite[Theorem~5.5]{CasasHerzogWachsmuth2015:1}
and this proof can be repaired by using
\cref{lem:weak_jop} above,
see also the arguments in the proof of \cref{lem: inequality j2 in case mu jp + F to 0}
below.
Thus,
\cite[Theorem~5.5]{CasasHerzogWachsmuth2015:1}
remains correct.

The next lemma will be used
to provide a lower bound for the second subderivative.
\begin{lemma}
	\label{lem: inequality j2 in case mu jp + F to 0}
	Let $\seq{v_k} \subset \LtwoOmegaT$ be a sequence which satisfies
	$v_k \weakly v \in \criticalcone{\bar u}$
	and $\mu j_2'(\bar u;v_k) + F'(\bar u)v_k \to 0$.
	Then
	\begin{equation*}
		\liminf_{k \to \infty} \int_0^T {{j_\Omega'(\bar u(t);v_k(t))}^2} \dt - j_2'(\bar u;v_k)^2
		\ge
		\int_0^T {{j_\Omega'(\bar u(t);v(t))}^2} \dt - j_2'(\bar u;v)^2
		.
	\end{equation*}
\end{lemma}
\begin{proof}
	It holds
	\begin{equation*}
		\mu \abs{j_2'(\bar u;v_k) - j_2'(\bar u;v)}
		\le
		\abs{\mu j_2'(\bar u;v_k) + F'(\bar u)v_k} +
		\abs{F'(\bar u) (v - v_k)} +
		\abs{-F'(\bar u) v - \mu j_2'(\bar u;v)}.
	\end{equation*}
	The first addend converges to zero by assumption,
	the second one by the weak convergence
	and the third one is zero as $v \in \criticalcone{\bar u}$.
	This implies $j_2'(\bar u;v_k) \to j_2'(\bar u;v)$.
	Hence,
	the subtrahend in the postulated inequality converges.
	For the minuend we use
	\begin{equation*}
		\int_0^T {{j_\Omega'(\bar u(t);v_k(t))}^2} \dt
		=
		\int_M {{j_\Omega'(\bar u(t);v_k(t))}^2} \dt
		+
		\int_A {{j_\Omega'(\bar u(t);v_k(t))}^2} \dt
		,
	\end{equation*}
	with $M$ as in \cref{lem:weak_jop}
	and $A := (0,T) \setminus M$.
	Combining \cref{lem:weak_jop} with the
	sequential weak lower semicontinuity
	of $\norm{\cdot}_{L^2(M)}^2$ yields
	\begin{equation*}
		\liminf_{k \to \infty} \int_M {{j_\Omega'(\bar u(t);v_k(t))}^2} \dt
		\ge
		\int_M {{j_\Omega'(\bar u(t);v(t))}^2} \dt
		.
	\end{equation*}
	Taking into account \eqref{jOmegap}, $j_\Omega'$ simplifies on the remaining set $A$
	and we get
	\begin{align*}
		\liminf_{k \to \infty} \int_A {{j_\Omega'(u(t);v_k(t))}^2} \dt
		&=
		\liminf_{k \to \infty} \int_A \parens*{ \int_\Omega \abs{v_k(x,t)} \dx }^2 \dt \\
		&\ge
		\int_A \parens*{ \int_\Omega \abs{v(x,t)} \dx }^2 \dt
		=
		\int_A {{j_\Omega'(u(t);v(t))}^2} \dt
		.
	\end{align*}
	The inner integral is continuous and convex as is also the square,
	leading to sequential weak lower semicontinuity.
	Adding both inequalities completes the proof.
\end{proof}

One could ask why we did it this way because $\int_0^T {{j_\Omega'(\bar u(t);v_k(t))}^2} \dt$ looks convex but this is not true.
The reason is that $j_\Omega'(\bar u(t);v_k(t))$ is convex w.r.t.\ $v_k$,
but squaring destroys convexity since $j_\Omega'(\bar u(t);v_k(t))$ can be negative.

The next lemma is similar,
but it provides an equality
whenever the sequence $v_k$ converges strongly.
\begin{lemma}
	\label{Integral jOmegaStrichQuadrat stetig}
	For $v_k \to v$ in $\LtwoOmegaT$
	we have
	$\int_0^T { j_\Omega'(\bar u;v_k)^2} \dt \to \int_0^T { j_\Omega'(\bar u;v)^2} \dt$.
\end{lemma}
\begin{proof}
	The estimates
	$\abs{j_\Omega'(\bar u;v_k) - j_\Omega'(\bar u;v)} \le \normLoneOmega{v_k - v}$
	and
	$\abs{j_\Omega'(\bar u;v)} \le \normLoneOmega{v}$
	follow easily by \eqref{jOmegap}.
	Using those estimates, we get
	\begin{align*}
			\MoveEqLeft
			\abs[\bigg]{\int_0^T { j_\Omega'(\bar u;v_k)^2 } \dt - \int_0^T { j_\Omega'(\bar u;v)^2 } \dt}
			\le
			\int_0^T { \abs*{j_\Omega'(\bar u;v_k)^2 - j_\Omega'(\bar u;v)^2} } \dt \\
			&\le
			\int_0^T { \abs*{ j_\Omega'(\bar u;v_k) - j_\Omega'(\bar u;v) } \abs*{ j_\Omega'(\bar u;v_k) }} \dt
			+ \int_0^T { \abs*{ j_\Omega'(\bar u;v_k) - j_\Omega'(\bar u;v) } \abs*{ j_\Omega'(\bar u;v) }} \dt \\
			&\le
			\sqrt{\int_0^T { \abs*{ j_\Omega'(\bar u;v_k) - j_\Omega'(\bar u;v) }^2} \dt}
			\parens[\Bigg]{
				\sqrt{\int_0^T { \abs*{ j_\Omega'(\bar u;v_k) }^2} \dt} + \sqrt{\int_0^T { \abs*{ j_\Omega'(\bar u;v) }^2} \dt}
			} \\
			&\le
			\sqrt{\int_0^T { \normLoneOmega{v_k - v}^2} \dt}
			\parens[\Bigg]{
			\sqrt{\int_0^T { \normLoneOmega{v_k}^2} \dt} + \sqrt{\int_0^T { \normLoneOmega{v}^2} \dt}
			} \\
			&=
			j_2(v_k - v)
			\parens*{ j_2(v_k) + j_2(v) }
			.
	\end{align*}
	The second inequality follows from the binomial formula,
	the third one by Hölder's inequality,
	the fourth one uses our previous estimates.
	Since $v_k \to v$, $j_2(v_k)$ is bounded
	and $j_2(v_k - v)$ converges to zero.
	This finishes the proof.
\end{proof}
The next lemma enables us to invoke \cref{Lemma bei j1 ueber F + mu j' = 0}
for the characterization of the critical cone $\criticalcone{\bar u}$.
\begin{lemma}
	\label{lem:equality_subdif_directional_2}
	Let $\lambda_{\bar u} \in \partial j_2(\bar u)$
	be given.
	For $v \in \LtwoOmegaT$,
	we have
	$\dual{\lambda_{\bar u}}{v} = j_2'(\bar u; v)$
	if and only if
	\begin{subequations}
		\label{subdifferential j2 charact}
		\begin{align}
			\label{subdifferential j2 charact_1}
			\bar u \ne 0: &\quad 
			\lambda_{\bar u}(x,t) \in \Sign(v(x,t)) \frac {\normLoneOmega{\bar u(t)}} {j_2(\bar u)} \text{ a.e.\ in } \set{\bar u = 0}
			\\
			\label{subdifferential j2 charact_2}
			\bar u = 0, v \ne 0: &\quad 
			\lambda_{\bar u}(x,t) \in \Sign(v(x,t)) \frac {\normLoneOmega{v(t)}} {j_2(v)} \text{ a.e.\ in } \Omega_T
			.
		\end{align}
	\end{subequations}
\end{lemma}
\begin{proof}
	We first consider the case $\bar u \ne 0$.
	From \eqref{eq:subdiff_j2} we already have
	$\lambda_{\bar u}(x,t) = s_{\bar u}(x,t) \norm{\bar u(t)}_{L^1(\Omega)} / j_2(\bar u)$
	with $s_{\bar u}(x,t) \in \Sign(\bar u(x,t))$
	for a.a.\ $(x,t) \in \Omega_T$.
	Further,
	\begin{equation*}
		j_2'(\bar u; v) - \dual{\lambda_{\bar u}}{v}
		=
		\frac1{j_2(\bar u)} \int_0^T \bracks*{
			j_\Omega'(\bar u(t); v(t)) - \int_\Omega s_{\bar u}(x,t) v(x,t) \d x
		} \norm{\bar u(t)}_{L^1(\Omega)}
		\d t.
	\end{equation*}
	Since the condition on $s_{\bar u}$ can be rewritten as
	$s_{\bar u}(t) \in \partial j_\Omega(\bar u(t))$,
	we can argue exactly as in \cref{lem:equality_subdif_directional_1}.

	It remains to consider $\bar u = 0$, $v \ne 0$.
	Due to $\bar u = 0$, we get $j_2'(\bar u; v) = j_2(v)$.
	Considering \eqref{eq:subdiff_j2} again,
	the condition in \eqref{subdifferential j2 charact_2}
	is equivalent to $\lambda_{\bar u} \in \partial j_2(v)$.
	Thus, it remains to show the equivalence of
	$\dual{\lambda_{\bar u}}{v} = j_2(v)$
	and
	$\lambda_{\bar u} \in \partial j_2(v)$.

	``$\Rightarrow$'':
	From \eqref{eq:subdiff_j2},
	we get $\norm{\lambda_{\bar u}}_{L^2(0,T;L^\infty(\Omega))} \le 1$.
	Thus, $j_2(w) \ge \dual{\lambda_{\bar u}}{w}$ for all $w \in \LtwoOmegaT$.
	Consequently, $j_2(w) - j_2(v) \ge \dual{\lambda_{\bar u}}{w - v}$ for all $w \in \LtwoOmegaT$.

	``$\Leftarrow$'':
	This follows from taking $w = 2 v$ and $w = 0$ in the subgradient inequality.
\end{proof}
The final lemma addresses the construction of a recovery sequence.
\begin{lemma}
	\label{lem:special_vk_for_j2}
	Assume $\bar u \ne 0$,
	$-F'(\bar u) \in \partial G_2(\bar u)$
	and let $v \in \criticalcone{\bar u}$ be given.
	For an arbitrary sequence $\seq{t_k} \subset \R^+$ with $t_k \searrow 0$,
	the sequence $\seq{v_k} \subset \LtwoOmegaT$ defined
	in
	\eqref{eq:special_vk_for_j1}
	satisfies
	\begin{subequations}
		\label{eq:special_vk_for_j2_properties}
		\begin{align}
			\label{eq:special_vk_for_j2:1}
			v_k &\to v \quad\text{in } \LtwoOmegaT, \\
			\label{eq:special_vk_for_j2:2}
			v_k &\in \criticalcone{\bar u}, \\
			\label{eq:special_vk_for_j2:3}
			\bar u + t_k v_k &\in \Uad, \\
			\label{eq:special_vk_for_j2:4}
			j_\Omega(\bar u + t_k v_k) - j_\Omega(\bar u) &= t_k j_\Omega'(\bar u;v_k) \quad\text{a.e.\ on } (0,T).
		\end{align}
	\end{subequations}
\end{lemma}
\begin{proof}
	As the sequence is the same as in \cref{lem:special_vk_for_j1},
	\eqref{eq:special_vk_for_j2:1} and \eqref{eq:special_vk_for_j2:3} have already been proven there,
	and
	\eqref{eq:special_vk_for_j2:4} can be shown analogously.

	It remains to verify \eqref{eq:special_vk_for_j2:2}.
	As in \cref{lem:special_vk_for_j1},
	we argue via \cref{Lemma bei j1 ueber F + mu j' = 0}
	and
	we analogously get that \eqref{Lemma j zu Zeigendes:1}
	is valid with $v$ replaced by $v_k$.
	It remains to show that
	$j_2'(\bar u; v) = \dual{\lambda_{\bar u}}{v}$
	implies
	$j_2'(\bar u; v_k) = \dual{\lambda_{\bar u}}{v_k}$.
	To this end,
	we can argue as in
	\cref{lem:special_vk_for_j1}
	by
	utilizing
	\eqref{subdifferential j2 charact_1}
	from \cref{lem:equality_subdif_directional_2}.
\end{proof}

Note that the case $\bar u = 0$
has been excluded in \cref{lem:special_vk_for_j2}.
The reason is that
condition
\eqref{subdifferential j2 charact_2}
is incompatible with pointwise changes
of $v$, see also \cref{ex:critical_cone_j2}
below.

Now, we are able to prove the main result of this section.

\begin{theorem}\label{G2''}
	We assume $-F'(\bar u) \in \partial G_2(\bar u)$.
	In case $\bar u \ne 0$,
	we have
	\begin{equation*}
		G_2''(\bar u, -F'(\bar u); v) =
		\frac \mu {j_2(\bar u)} \parens[\bigg]{ \int_0^T {j_\Omega'(\bar u(t);v(t))^2} \dt - j_2'(\bar u;v)^2 }
	\end{equation*}
	for all $v \in \criticalcone{\bar u}$
	and $G_2$ is strongly twice epi-differentiable at $\bar u$ for $-F'(\bar u)$.

	In case $\bar u = 0$,
	we have
	\begin{equation*}
		G_2''(\bar u, -F'(\bar u); v) \ge 0
		\quad\forall v \in \criticalcone{\bar u}
		\qquad\text{and}\qquad
		G_2''(\bar u, -F'(\bar u); v) = \infty
		\quad\forall v \in \LtwoOmegaT \setminus \criticalcone{\bar u}
		.
	\end{equation*}
\end{theorem}
\begin{proof}
	We first consider $\bar u \ne 0$.
	Let $v \in \criticalcone{\bar u}$ be given.
	The first step is to show that the above right-hand side
	is a lower bound for the expression
	\begin{equation*}
		L :=
		\liminf_{k \to \infty}
		\frac{G_2(\bar u + t_k v_k) - G_2(\bar u) - t_k \dual{-F'(\bar u)}{v_k}}{t_k ^2 / 2} \\
	\end{equation*}
	for every pair of sequences $\seq{t_k} \subset \R^+$ and $\seq{v_k} \subset \LtwoOmegaT$
	with $t_k \searrow 0$ and $v_k \weakly v$.
	It is clear that we only have to consider sequences with
	$\bar u + t_k v_k \in \Uad$.
	As $\seq{t_k}$ is a zero sequence and $\seq {v_k}$ converges weakly,
	\eqref{taylor estimate for j2} holds $k$ large enough.
	We use the abbreviation
	\begin{equation*}
		\Theta(\bar u,v_k)
		:=
		\frac 1 {j_2(\bar u)} \parens[\bigg]{ \int_0^T {{j_\Omega'(\bar u(t);v_k(t))}^2} \dt - j_2'(\bar u;v_k)^2 }
		.
	\end{equation*}
	Note that $\Theta(\bar u, v_k) \ge 0$
	due to Hölder's inequality.
	With \cref{lem:lower_taylor_expansion}, we get
	\begin{align*}
		L
		&=
		\liminf_{k \to \infty}
		\frac 2 {t_k^2}
		\parens[\Big]{
			\mu \bracks*{ j_2(\bar u + t_k v_k) - j_2(\bar u) } + t_k \dual{F'(\bar u)}{v_k}
		} \\
		&\ge
		\liminf_{k \to \infty}
		\frac 2 {t_k^2}
		\parens*{
			\mu \bracks*{
				t_k j_2'(\bar u;v_k) + \frac {t_k^2} 2 \Theta(\bar u, v_k) - \frac {C t_k^3 \Ltwonorm{v_k}^3} {j_2(\bar u)^2}
			}
			+ t_k \dual{F'(\bar u)}{v_k}
		}.
	\end{align*}
	As $\Ltwonorm{v_k}$ is bounded and $j_2(\bar u) > 0$ holds, the cubic term in brackets vanishes
	as $k \to \infty$.
	This yields
	\begin{equation*}
		L \ge
		\liminf_{k \to \infty}
		\parens*{
			\frac 2 {t_k} \parens[\Big]{
				\mu j_2'(\bar u;v_k) + \dual{F'(\bar u)}{v_k}
			}
			+ \mu \Theta(\bar u, v_k)
		}.
	\end{equation*}
	Note that
	$-F'(\bar u) \in \partial G_2(\bar u)$
	and $\bar u + t_k v_k \in \Uad$
	yield
	$\mu j_2'(\bar u;v_k) + \dual{F'(\bar u)}{v_k} \ge 0$.
	Now, we distinguish two cases.

	\emph{Case 1}:
	$\liminf_{k \to \infty} \mu j_2'(\bar u;v_k) + \dual{F'(\bar u)}{v_k} > 0$.
	Because of the factor $2/t_k$ we get $L = \infty$.
	The desired inequality is fulfilled.

	\emph{Case 2}:
	$\liminf_{k \to \infty} \mu j_2'(\bar u;v_k) + \dual{F'(\bar u)}{v_k} = 0$.
	We choose subsequences of $\seq{t_k}$ and $\seq{v_k}$ (without relabeling)
	which realize the limit inferior,
	i.e.,
	$\mu j_2'(\bar u; v_k) + \dual{F'(\bar u)}{v_k} \to 0$.
	In this situation,
	\cref{lem: inequality j2 in case mu jp + F to 0}
	can be applied and yields
	the desired
	\begin{equation*}
		L
		\ge
		\frac \mu {j_2(\bar u)}
		\parens[\bigg]{
			\int_0^T {{j_\Omega'(\bar u(t);v(t))}^2} \dt - j_2'(\bar u;v)^2
		}.
	\end{equation*}

	We will now show that this lower bound is realized
	for an arbitrary sequence $\seq{t_k} \subset \R^+$ with $t_k \searrow 0$
	if $\seq{v_k}$ is chosen as in \cref{lem:special_vk_for_j2}.
	For the purpose of shortening, let $S(\bar u, k) := j_2(\bar u + t_k v_k) + j_2(\bar u)$.
	We note that
	$2/S(\bar u, k) \to 1/ j_2(\bar u)$. We get
	\begingroup
	\allowdisplaybreaks
	\begin{align*}
		\MoveEqLeft
		\lim_{k \to \infty}
		\frac{G_2(\bar u + t_k v_k) - G_2(\bar u) - t_k \dual{-F'(\bar u)}{v_k}}{t_k ^2 / 2}
		\\
		&=
		\lim_{k \to \infty}
		\frac 2 {t_k^2} \parens*{ \mu \bracks*{ j_2(\bar u + t_k v_k) - j_2(\bar u) } + t_k F'(\bar u)v_k }
		&&\mathllap{\text{(by \eqref{eq:special_vk_for_j2:3})}}
		\\
		&=
		\lim_{k \to \infty}
		\frac
		{2 \parens*{ \mu \bracks*{ j_2(\bar u + t_k v_k)^2 - j_2(\bar u)^2 }
		+ t_k S(\bar u, k) F'(\bar u)v_k }}
		{t_k^2 S(\bar u, k)}
		\\
		&=
		\frac 1 { j_2(\bar u)}
		\lim_{k \to \infty}
		\frac
		{
			\mu \int_0^T {j_\Omega(\bar u + t_k v_k)^2 - j_\Omega(\bar u)^2} \dt
			+ t_k S(\bar u, k) F'(\bar u)v_k
		}
		{t_k^2}
		&\mspace{90mu}&\mathllap{\text{(by \eqref{eq: other form of j2})}}
		\\
		&=
		\frac 1 { j_2(\bar u)}
		\lim_{k \to \infty}
		\frac
		{
			\mu \int_0^T { t_k j_\Omega'(\bar u;v_k)^2 + 2 j_\Omega(\bar u) j_\Omega'(\bar u;v_k) } \dt
			+
			S(\bar u, k) F'(\bar u)v_k
		}
		{t_k}
		&&\mathllap{\text{(by \eqref{eq:special_vk_for_j2:4})}}
		\\
		&=
		\frac \mu { j_2(\bar u)}
		\lim_{k \to \infty}
		\frac
		{
			\int_0^T { t_k j_\Omega'(\bar u;v_k)^2 } + 2 j_\Omega(\bar u) j_\Omega'(\bar u;v_k) \dt
			-
			S(\bar u, k) j_2'(\bar u;v_k)
		}
		{t_k}
		&&\mathllap{\text{(by \eqref{eq:special_vk_for_j2:2})}}
		\\
		&=
		\frac \mu { j_2(\bar u)}
		\lim_{k \to \infty} \frac 1 {t_k}
		\parens[\bigg]{
			\int_0^T { t_k j_\Omega'(\bar u;v_k)^2 } \dt
			+
			\bracks*{
				2 j_2(\bar u) - S(\bar u, k)
			}
			j_2'(\bar u;v_k)
		}
		&&\mathllap{\text{(by \eqref{eq:directional_derivative_j2})}}
		\\
		&=
		\frac \mu { j_2(\bar u)}
		\lim_{k \to \infty}
		\parens[\bigg]{
			\int_0^T { j_\Omega'(\bar u;v_k)^2 } \dt
			- \frac {j_2(\bar u + t_k v_k) - j_2(\bar u)} {t_k} j_2'(\bar u;v_k)
		}.	
	\end{align*}
	\endgroup
	The expression $j_2'(\bar u;\cdot)$ is continuous and
	$\frac {j_2(\bar u + t_k v_k) - j_2(\bar u)} {t_k} \to j_2'(\bar u;v)$ holds
	(cf.\ \cref{lem:sumrule_subderivative_take_2} as $j_2$ is convex and Lipschitz continuous and therefore Hadamard differentiable).
	Together with
	\Cref{Integral jOmegaStrichQuadrat stetig}, this yields the claim.

	In case $\bar u = 0$,
	the assertion directly follows from
	\cref{lem:Gpp_convex} and \eqref{eq:heile_welt_on_KK}.
\end{proof}

The next example shows that the situation $\bar u = 0$
is surprisingly difficult, even in case $F'(\bar u) \in L^\infty(\Omega_T)$.
\begin{example}
	\label{ex:critical_cone_j2}
	We use the setting $T = 1$, $\Omega = (0,1)$, i.e., $\OmegaT = (0,1)^2$.
	Further, for some $\rho \in (0,1)$ we set
	\begin{equation*}
		D :=
		\set{
			(x,t) \in \OmegaT
			\given
			0 < x < t^\rho < 1
		}.
	\end{equation*}
	Next, we fix $\bar u = 0$
	and we assume that the smooth part of the objective satisfies
	$F'(\bar u) \equiv -1$ on $D$ while $\abs{F'(\bar u)} < 1$ on $\OmegaT \setminus D$.
	Finally, we set $\alpha := -1$, $\beta := 1$ and $\mu = 1$.

	First, we show that the critical cone is nonempty.
	We define
	the measurable function
	$v \colon \Omega_T \to \R$ via
	\begin{equation*}
		v(x,t) :=
		\begin{cases}
			t^{-\rho} & \text{if } (x,t) \in D \\
			0 & \text{else}.
		\end{cases}
	\end{equation*}
		Due to
		\begin{equation*}
			\Ltwonorm{v}^2
			=
			\int_0^T \int_0^{t^\rho} t^{-2\rho} \dx \dt
			=
			\int_0^T t^{-\rho} \dt
			=
			\frac 1 {1-\rho}
			,
		\end{equation*}
	we have $v \in \LtwoOmegaT$.
	Now, we can check
	\begin{equation*}
		j_2'(\bar u; v)
		=
		j_2(v)
		=
		\parens[\bigg]{
			\int_0^T \norm{v(t)}_{L^1(\Omega)}^2 \, \dt
		}^{1/2}
		=
		1
	\end{equation*}
	and
	\begin{equation*}
		-F'(\bar u) v
		=
		\int_D v \d(x,t)
		=
		\int_0^T \int_0^{t^\rho} \frac1 {t^\rho} \d x \d t
		=
		1.
	\end{equation*}
	Thus, $F'(\bar u) v + j_2'(\bar u; v) = 0$
	and we trivially have $v \in \tangentcone{\bar u} = L^2(\Omega_T)$.
	Hence, $v \in \criticalcone{\bar u}$.

	Next, we show that every $\tilde v \in \criticalcone{\bar u} \setminus \set{0}$
	is unbounded, i.e., $\tilde v \not\in L^\infty(\OmegaT)$.
	Indeed, for an arbitrary $\tilde v \in \criticalcone{\bar u} \setminus \set{0}$ we have
	\begin{equation*}
		j_2(\tilde v)
		=
		j_2'(\bar u; \tilde v)
		=
		-F'(\bar u) \tilde v
		\le
		\int_0^T \norm{\tilde v(t)}_{L^1(\Omega)} \d t
		\le
		\norm{ \norm{\tilde v(t)}_{L^1(\Omega)} }_{L^2(0,T)}
		=
		j_2(\tilde v).
	\end{equation*}
	Hence, both inequalities are actually equalities.
	This implies $\tilde v \ge 0$ a.e.\ on $D$, $\tilde v = 0$ a.e.\ on $\Omega_T \setminus D$
	and that $\norm{\tilde v(t)}_{L^1(\Omega)}$ is a constant, say, $c > 0$.
	For a.e.\ $t \in (0,1)$, we get $c = \int_\Omega \tilde v(x,t) \d x = \int_0^{t^\rho} \tilde v(x,t) \d x$
	and this shows
	$\measure*{1}{\set{x \in \Omega \given \tilde v(x,t) \ge c/t^\rho }} > 0$,
	where $\lambda^d$ denotes the $d$-dimensional Lebesgue measure.
	Consequently, Fubini implies that for any $\tau \in (0,1)$,
	we have
	\begin{align*}
		\measure*{2}{\set{(x,t) \in \Omega_T \given \tilde v(x,t) \ge c/\tau }}
		&=
		\int_0^T \measure*{1}{\set{ x \in \Omega \given \tilde v(x,t) \ge c/\tau }} \d t
		\\&\ge                                                                          
		\int_0^{\tau^{1/\rho}} \measure*{1}{\set{ x \in \Omega \given \tilde v(x,t) \ge c/t^\rho }} \d t
		>
		0.
	\end{align*}
	This shows that $\tilde v \not\in L^\infty(\Omega_T)$.
	In particular,
	for any $\tilde v \in \criticalcone{\bar u} \setminus \set{0}$ and any $t > 0$,
	we have $\bar u + t \tilde v \not\in \Uad$.
	Further, this shows that the assertion of \cref{lem:special_vk_for_j2}
		is not valid in case $\bar u = 0$, even if $F'(\bar u) \in \LinftyOmegaT$.

	Finally, let a sequence $\seq{t_k} \subset \R^+$ with $t_k \searrow 0$ be arbitrary
	and we consider the approximation
	$v_k := P_{[-1/t_k, 1/t_k]}(v)$ of the above $v$.
	It is clear that $v_k \to v$ in $\LtwoOmegaT$
	and $\bar u + t_k v_k \in \Uad$.
	A simple calculation shows (for $t_k$ small enough)
	\begin{align*}
		j_2(v_k)
		&=
		\parens[\bigg]{
			\int_0^{t_k^{1/\rho}} \norm{v_k(t)}_{L^1(\Omega)}^2 \d t
			+
			\int_{t_k^{1/\rho}}^1 \norm{v_k(t)}_{L^1(\Omega)}^2 \d t
		}^{\frac12}
		=
		\parens[\bigg]{
			\int_0^{t_k^{1/\rho}} \frac{t^{2\rho}}{t_k^2}\d t
			+
			\int_{t_k^{1/\rho}}^1 1 \d t
		}^{\frac12}
		\\&
		=
		\parens*{
			\frac{t_k^{1/\rho}}{2\rho+1} + 1 - t_k^{1/\rho}
		}^{1/2}
		=
		\parens*{ 1 - (2\rho/(2\rho+1)) t_k^{1/\rho}}^{1/2}
	\end{align*}
	and
	\begin{equation*}
		-F'(\bar u) v_k
		=
		\int_0^1 \int_\Omega v_k \d x \d t
		=
		\int_0^{t_k^{1/\rho}} t^\rho/t_k \d t
		+
		\int_{t_k^{1/\rho}}^1 1 \d x \d t
		=
		1 - \frac {\rho} {\rho+1} t_k^{1/\rho}
		.
	\end{equation*}
	Thus,
	the curvature term is
	\begin{align*}
		\lim_{k \to \infty} \frac{ G_2(\bar u + t_k v_k) - G_2(\bar u) + t_k F'(\bar u) v_k}{t_k^2/2}
		&=
		\lim_{k \to \infty} \frac{ \parens*{ 1 - \frac {2\rho} {2\rho+1} t_k^{1/\rho}}^{1/2} - \left( 1 - \frac {\rho} {\rho+1} t_k^{1/\rho} \right)}{t_k/2}
		=
		0.
	\end{align*}
	This shows that $G_2''(\bar u, -F'(\bar u); v) = 0$
		and that $\seq{v_k}$ serves as
		a recovery sequence.
		However, it is not clear
		whether a similar approach works for all $\tilde v \in \criticalcone{\bar u}$
		and whether $G_2''(\bar u, -F'(\bar u); \tilde v) = 0$ for all $\tilde v \in \criticalcone{\bar u}$.
\end{example}

\subsection{Second subderivative of \texorpdfstring{$j_3$}{j3}}

As in \cite[p. 273,292]{CasasHerzogWachsmuth2015:1},
we define the sets
\begin{subequations}
	\begin{align}
		\Omega_{\bar u} &:= \set{ x \in \Omega \given \normLtwotime{\bar u(x)} \ne 0 },\\
		\Omega_{\bar u}^0 &:= \set{ x \in \Omega \given \normLtwotime{\bar u(x)} = 0 } = \Omega \setminus \Omega_{\bar u},\\
		\label{eq:omega_sigma}
		\Omega_\sigma &:= \set{ x \in \Omega \given \normLtwotime{\bar u(x)} \ge \sigma },
		\qquad\forall \sigma > 0.
	\end{align}
\end{subequations}
From
\cite[Proposition~3.8]{CasasHerzogWachsmuth2015:1},
we recall
the directional derivative of $j_3$
\begin{equation}
	\label{eq:directional_derivative_j3}
	j_3'(\bar u; v) =
	\int_{\Omega_{\bar u}^0} {\normLtwotime{v(x)}} \dx +
	\int_{\Omega_{\bar u}} {\frac 1 {\normLtwotime{\bar u(x)}} \int_0^T {\bar u v} \dt } \dx
\end{equation}
and $\lambda \in \partial j_3(\bar u)$ is equivalent to
$\lambda \in L^\infty(\Omega; L^2(0,T))$ and
\begin{subequations}
	\label{lem:subdifferential_of_j3}
	\begin{align}
		\label{lem:subdifferential_of_j3_part_a}
		\normLtwotime{\lambda(x)} &\le 1
		&&
		\text{f.a.a.\ } x \in \Omega_{\bar u}^0
		\\
		\label{lem:subdifferential_of_j3_part_b}
		\lambda(x,t) &= \frac{\bar u(x,t)}{\normLtwotime{\bar u(x)}}
		&&
		\text{f.a.a.\ } x \in \Omega_{\bar u} \text{ and } t \in (0,T).
	\end{align}
\end{subequations}

\begin{lemma}\label{lem:j3_first_part_of_proof:swlsc}
	For any measurable $M \subset \Omega_{\bar u}$,
	the mapping
	$q_M \colon \LtwoOmegaT \to \R$
	given by
	\begin{equation*}
		q_M(v) :=
		\int_M {\frac 1 {\normLtwotime{\bar u(x)}} \bracks*{ \int_0^T {v^2(x,t)} \dt - \parens*{ \int_0^T {\frac{\bar u(x,t) v(x,t)}{\normLtwotime{\bar u(x)}}} \dt }^2 }} \dx
	\end{equation*}
	is convex, lower semicontinuous,
	and therefore sequentially weakly lower semicontinuous.
\end{lemma}
\begin{proof}
	For fixed $\sigma > 0$,
	we define $M_\sigma := M \cap \Omega_\sigma$ and
	$b_\sigma : \LtwoOmegaT^2 \to \R$,
	\begin{equation*}
		b_\sigma(v,w) :=
		\int_{M_\sigma} \frac 1 {\normLtwotime{\bar u(x)}} \bracks*{ \int_0^T {vw} \dt - \frac{1}{\normLtwotime{\bar u(x)}^2}\parens*{ \int_0^T \bar u v \dt } \parens*{ \int_0^T \bar u w \dt } } \dx.
	\end{equation*}
	This is a symmetric and real-valued bilinear form.
	From Hölder's inequality
	we get
	\begin{equation}\label{eq:integrand_of_j3_nonnegative}
		\frac 1 {\normLtwotime{\bar u(x)}}
		\bracks*{
			\int_0^T {v^2(x,t)} \dt
			- \parens*{ \int_0^T {\frac{\bar u(x,t) v(x,t)}{\normLtwotime{\bar u(x)}}} \dt }^2
		}
		\ge 0
	\end{equation}
	for all $x \in \Omega_{\bar u}$.
	Hence,
	$b_\sigma(v,v) \ge 0$
	and, therefore,
	$v \mapsto b_\sigma(v,v)$
	is convex.
	For the continuity of $b_\sigma$,
	we note
	\begin{equation*}
		\abs{b_\sigma(v,v)}
		\le
		\int_{M_\sigma}
		\frac {\normLtwotime{v}^2} {\normLtwotime{\bar u(x)}}
		\dx
		\le
		\frac{\Ltwonorm{v}^2}{\sigma}
		.
	\end{equation*}
	Together with the symmetry, we get that $b_\sigma$ is bounded, hence continuous.

	The monotone convergence theorem
	yields $q(v) = \lim_{\sigma \searrow 0} b_\sigma(v,v) = \sup_{\sigma > 0} b_\sigma(v,v)$.
	Since
	the supremum of convex and lower semicontinuous functions is again convex and lower semicontinuous,
	this establishes the claim.
\end{proof}

The next lemma follows
since $L^2(0,T)$ is a Hilbert space.
\begin{lemma}
	\label{lem:psi_and_its_derivatives}
	We define the function $\Psi: L^2(0, T) \to \R$ by $\Psi(f) := \normLtwotime{f}$.
	For every $f \ne 0$ and $g \in L^2(0,T)$, we have
	\begin{subequations}
		\begin{align}
			\Psi'(f) g &= \frac 1 {\normLtwotime{f}} \int_0^T {fg} \dt
			,
			\\
			\Psi''(f) g^2 &= \frac 1 {\normLtwotime{f}} \braces*{ \int_0^T {g^2} \dt - \frac 1 {\normLtwotime{f}^2} \parens*{ \int_0^T {fg} \dt }^2 }
			,
			\\
			\Psi'''(f)g^3 &= \frac 3 {\normLtwotime{f}^3} \braces*{ \frac 1 {\normLtwotime{f}^2} \parens*{ \int_0^T {fg} \dt }^3 - \parens*{ \int_0^T {g^2} \dt } \parens*{ \int_0^T {fg} \dt } }
			.
		\end{align}
		Furthermore,
		\begin{equation}
			\label{eq:psi_and_its_derivatives:4}
			\abs*{ \Psi'''(f)g^3 } \le  \frac {6 \normLtwotime{g}^3} {\normLtwotime{f}^2}.
		\end{equation}
	\end{subequations}
\end{lemma}

The next lemma will be used
to show a lower bound for the second subderivative of $j_3$.

\begin{lemma}\label{j3 Beweisteil 1 Abschaetzung}
	We assume $\bar u \ne 0$ and
	let sequences $\seq{t_k} \subset \R^+$ and $\seq{v_k} \subset \LtwoOmegaT$
	be given such that
	$t_k \searrow 0$ and $v_k \weakly v$.
	Then, it holds
	\begin{equation}
		\label{eq:j3_first_part_of_proof:estimate_for_j3pp-term}
		\begin{aligned}
			\MoveEqLeft
			\liminf_{k \to \infty}
			\frac{j_3(\bar u + t_k v_k) - j_3(\bar u) - t_k j_3'(\bar u; v_k)}{t_k^2 / 2} \\
			&\ge
			\int_{\Omega_{\bar u}} 
			{\frac 1 {\normLtwotime{\bar u(x)}} \braces*{ \int_0^T {v(x,t)^2} \dt - \frac {\parens*{ \int_0^T {\bar u(x,t) v(x,t)} \dt }^2} {\normLtwotime{\bar u(x)}^2} }}
			\dx
			.
		\end{aligned}
	\end{equation}
\end{lemma}
\begin{proof}
	First we extract subsequences that realize the limit inferior,
	afterwards we extract subsequences (again without relabeling) such that
	\begin{equation}
		\label{eq:bound_on_t_k}
		\forall k \in \N: t_k \le \frac 1 {k^4}.
	\end{equation}
	For every $N \in \N$ we define the set $M_N$ and the functional $j_{3,N}$ via
	\begin{align*}
		M_N &:=
		\set*{
			x \in \Omega_{1/N} \given
			\forall k \ge N :
			\normLtwotime{v_k(x)} \le
			t_k^{-1/4}
		}
		,
		\\
		j_{3,N}(u) &:=
		\int_{M_N} \normLtwotime{u(x, \cdot)} \dx
		.
	\end{align*}
	Our first goal is to show that the analogue of \eqref{eq:j3_first_part_of_proof:estimate_for_j3pp-term}
	holds for the functional $j_{3,N}$
	for fixed $N \in \N$ with $N \ge 2$.
	For all $x \in M_N$, $k \ge N$ and $\theta \in [0,1]$
	we get
	with \eqref{eq:omega_sigma} and \eqref{eq:bound_on_t_k}
	\begin{align*}
			\normLtwotime{\bar u(x) + \theta t_k v_k(x)}
			&\ge
			\normLtwotime{\bar u(x)} - \theta \normLtwotime{t_k v_k(x)}
			\\&
			\ge
			\frac1N
			- t_k \norm{v_k(x)}_{L^2(0,T)}
			\ge
			\frac1N - t_k^{3/4}
			\ge
			\frac1N - \frac{1}{k^3}
			\ge
			\frac1{2 N}
	\end{align*}
	and therefore
	\begin{equation*}
		0
		\le
		\frac {t_k \normLtwotime{v_k(x)}^3} {\normLtwotime{\bar u(x) + \theta t_k v_k(x)}^2}
		\le
		4 N^2 t_k^{1/4}.
	\end{equation*}
	Using
	$\Psi$ from \Cref{lem:psi_and_its_derivatives}
	we perform a Taylor expansion
	and obtain $\theta \in [0,1]$ (depending on $k \ge N$ and $x \in M_N$)
	such that
	\begin{align*}
		\MoveEqLeft
		\Psi(\bar u(x) + t_k v_k(x))
		-
		\Psi(\bar u(x))
		-
		t_k \Psi'(\bar u(x)) v_k(x)
		\\&
		=
		\frac{t_k^2}{2} \Psi''(\bar u(x)) v_k(x)^2
		+
		\frac{t_k^3}{6} \Psi'''(\bar u(x) + \theta t_k v_k(x)) v_k(x)^3
		.
	\end{align*}
	Together with \eqref{eq:psi_and_its_derivatives:4}
	and the above estimate,
	we get
	\begin{equation*}
		\Psi(\bar u(x) + t_k v_k(x))
		-
		\Psi(\bar u(x))
		-
		t_k \Psi'(\bar u(x)) v_k(x)
		\ge
		\frac{t_k^2}{2} \Psi''(\bar u(x)) v_k(x)^2
		-
		4 t_k^2 N^2 t_k^{1/4}
		.
	\end{equation*}
	Since this estimate holds for all $x \in M_N$,
	we can integrate and obtain
	\begin{equation*}
		j_{3,N}(\bar u + t_k v_k)
		-
		j_{3,N}(\bar u)
		-
		t_k j_{3,N}'(\bar u) v_k
		\ge
		\frac{t_k^2}{2} q_{M_N}(v_k)
		-
		4 t_k^2 N^2 t_k^{1/4}
		\int_{M_N} 1 \dx
	\end{equation*}
	with $q_{M_N}$ from \cref{lem:j3_first_part_of_proof:swlsc}.
	Now, we divide by $t_k^2/2$
	and, using \cref{lem:j3_first_part_of_proof:swlsc},
	we pass to the limit $k \to \infty$
	to obtain
	\begin{equation*}
		\liminf_{k \to \infty}\frac{
			j_{3,N}(\bar u + t_k v_k)
			-
			j_{3,N}(\bar u)
			-
			t_k j_{3,N}'(\bar u) v_k
		}{t_k^2 / 2}
		\ge
		q_{M_N}(v).
	\end{equation*}
	Since $j_3 - j_{3,N}$
	is a convex function,
	we get
	\begin{align*}
		\MoveEqLeft
		\liminf_{k \to \infty}\frac{
			j_{3}(\bar u + t_k v_k)
			-
			j_{3}(\bar u)
			-
			t_k j_{3}'(\bar u) v_k
		}{t_k^2 / 2}
		\\
		&
		\ge
		\liminf_{k \to \infty}\frac{
			j_{3,N}(\bar u + t_k v_k)
			-
			j_{3,N}(\bar u)
			-
			t_k j_{3,N}'(\bar u) v_k
		}{t_k^2 / 2}
		\ge
		q_{M_N}(v)
		.
	\end{align*}
	It remains to pass to the limit $N \to \infty$.
	Note that the set $M_N$ is increasing in $N$.
	Moreover,
	for the Lebesgue measure of $\Omega_{\bar u} \setminus M_N$ we get
	\begin{align*}
		\measure{d}{\Omega_{\bar u} \setminus M_N}
		&=
		\measure*{d}{\Omega_{\bar u} \setminus \Omega_{1/N}}
		+
		\measure*{d}{\Omega_{1/N} \setminus M_N}
		\\
		&=
		\measure*{d}{\Omega_{\bar u} \setminus \Omega_{1/N}}
		+
		\measure*{d}{\set*{x \in \Omega_{1/N} \given \exists k \ge N:
		\normLtwotime{v_k(x)} > t_k^{-1/4}}}
		\\
		&\le
		\measure*{d}{\Omega_{\bar u} \setminus \Omega_{1/N}}
		+
		\sum_{k \ge N}\measure*{d}{\set*{x \in \Omega \given
		\normLtwotime{v_k(x)} > t_k^{-1/4}}}
		.
	\end{align*}
	Next, we use Chebyshev's inequality
	and \eqref{eq:bound_on_t_k}
	to get
	\begin{align*}
		\measure{d}{\Omega_{\bar u} \setminus M_N}
		&\le
		\measure*{d}{\Omega_{\bar u} \setminus \Omega_{1/N}}
		+
		\sum_{k \ge N}
		t_k^{1/2}
		\int_{\Omega} \normLtwotime{v_k(x)}^2 \dx
		\\
		&\le
		\measure*{d}{\Omega_{\bar u} \setminus \Omega_{1/N}}
		+
		\sum_{k \ge N}
		k^{-2}
		\Ltwonorm{v_k}^2
		.
	\end{align*}
	The first addend trivially vanishes for $N \to \infty$.
	Furthermore, $\Ltwonorm{v_k}$ is bounded due to weak convergence and
	the series $\sum_{n \ge 1} n^{-2}$ converges absolutely.
	Therefore, we get the convergence $\measure{d}{\Omega_{\bar u} \setminus M_N} \to 0$ for $N \to \infty$.
	
	In order to finish the proof,
	we write
	\begin{align*}
		\MoveEqLeft
		\liminf_{k \to \infty}
		\frac{j_3(\bar u + t_k v_k) - j_3(\bar u) - t_k j_3'(\bar u; v_k)}{t_k^2 / 2}
		\\
		&\ge
		q_{M_N}(v)
		=
		\int_{\Omega_{\bar u}}
		{\frac {\characteristicfunction{M_N}(x)} {\normLtwotime{\bar u(x)}} \braces*{ \int_0^T {v^2} \dt - \frac {\parens*{ \int_0^T {\bar u v} \dt }^2} {\normLtwotime{\bar u(x)}^2} }}
		\dx
		.
	\end{align*}
	It is easy to see by \eqref{eq:integrand_of_j3_nonnegative}
	that the integrand is nonnegative for every $N \in \N$.
	Since $M_N$ is increasing,
	the sequence of integrands
	converges monotonely
	towards the integrand in \eqref{eq:j3_first_part_of_proof:estimate_for_j3pp-term}.
	An application of the monotone convergence theorem finishes the proof.
\end{proof}

The next lemma prepares the application of \cref{Lemma bei j1 ueber F + mu j' = 0},
in particular, it characterizes \eqref{Lemma j zu Zeigendes:2}
for $j = j_3$.
For convenience, we recall that the directional derivative of $j_3$
was given in \eqref{eq:directional_derivative_j3}.
\begin{lemma}
	\label{lem:equality_subdif_directional_3}
	Let $\lambda_{\bar u} \in \partial j_3(\bar u)$
	be given.
	For $v \in \LtwoOmegaT$,
	we have
	$\dual{\lambda_{\bar u}}{v} = j_3'(\bar u; v)$
	if and only if
	\begin{align}
		\label{subdifferential j3 charact_1}
		\lambda_{\bar u}(x,t) = \frac{v(x,t)}{\normLtwotime{v(x)}}
		\quad
		\text{f.a.a.\ } x \in \Omega_{\bar u}^0 \text{ with } \normLtwotime{v(x)} \ne 0
		.
	\end{align}
\end{lemma}
\begin{proof}
	The implication ``$\Leftarrow$'' is clear.

	Let $\dual{\lambda_{\bar u}}{v} = j_3'(\bar u; v)$ be fulfilled.
	With \eqref{lem:subdifferential_of_j3} to get
	\begin{align*}
		0 = 
		j_3'(\bar u; v) - \dual{\lambda_{\bar u}}{v}
		&=
		\int_{\Omega_{\bar u}^0} {\normLtwotime{v(x)}} - \int_0^T \lambda_{\bar u} v \dt \dx
		\\
		&\ge
		\int_{\Omega_{\bar u}^0} {\normLtwotime{v(x)} (1 - \normLtwotime{\lambda_{\bar u}(x)})} \dx
		\ge
		0
		.
	\end{align*}
	Hence,
	$\normLtwotime{v(x)}  \normLtwotime{\lambda_{\bar u}(x)} = \int_0^T \lambda_{\bar u} v \dt$
	and
	$\normLtwotime{v(x)} (1 - \normLtwotime{\lambda_{\bar u}(x)}) = 0$
	for a.a.\ $x \in \Omega_{\bar u}^0$.
	This yields \eqref{subdifferential j3 charact_1}.
\end{proof}

For $j_3$,
we cannot use the same construction \eqref{eq:special_vk_for_j1}
as for $j_1$ and $j_2$,
since this would lead to problems
on the set $\Omega_{\bar u}^0$, cf.\ \cref{ex:critical_cone_j2}.
In order to get $v_k \in C_{\bar u}$,
we have to modify the construction.
We follow
\cite[Theorem~4.3, Case III]{CasasHerzogWachsmuth2015:1}.
\begin{lemma}
	\label{lem:special_vk_for_j3}
	We assume $-F'(\bar u) \in \partial G(\bar u) \cap \LinftyOmegaT$.
	Let $\seq{t_k} \subset \R^+$ be an arbitrary sequence with $t_k \searrow 0$
	and $v \in \criticalcone{\bar u}$.
	We define the sequence 
	$\seq{v_k} \subset \LtwoOmegaT$
	on the set $\Omega_{\bar u} \times (0,T)$
	via
	\begin{align*}
		v_k := 
		\begin{cases}
			0
			&\text{if } \bar u \in (\alpha, \alpha +\sqrt{t_k}) \cup (\beta - \sqrt{t_k}, \beta) \cup (-\sqrt{t_k}, 0) \cup (0, \sqrt{t_k}),
			\\
			0
			&\text{if } \normLtwotime{\bar u(x)} < \sqrt{t_k},
			\\
			P_{t_k}(v)
			&\text{otherwise,}
		\end{cases}
	\end{align*}
	where $P_{t_k} \colon \R \to \R$ denotes the projection onto the interval
		$\bracks*{ -\frac 1 {\sqrt{t_k}}, \frac 1 {\sqrt{t_k}} }$,
	and on $\Omega_{\bar u}^0 \times (0,T)$
	we define
	\begin{align*}
		v_k :=
		\begin{cases}
			0
			&\text{if } \normLtwotime{v(x)} > \frac 1 {\sqrt{t_k}}
			\\
			v
			&\text{otherwise}.
		\end{cases}
	\end{align*}
	Then, this sequence satisfies
	(for $k$ large enough)
	\begin{subequations}
		\label{eq:special_vk_for_j3_properties}
		\begin{align}
			\label{eq:special_vk_for_j3:1}
			v_k &\to v \quad\text{in } \LtwoOmegaT,
			\\
			\label{eq:special_vk_for_j3:2}
			v_k &\in \criticalcone{\bar u},
			\\
			\label{eq:special_vk_for_j3:3}
			\bar u + t_k v_k &\in \Uad,
			\\
			\label{eq:special_vk_for_j3:4}
			\normLtwotime{\bar u(x) + t_k v_k(x)} - \normLtwotime{\bar u(x)}
			&=
			t_k \frac{\int_0^T {2 \bar u v_k + t_k v_k^2} \dt}{K_k(x)}
			\quad\text{for } x \in \Omega_{\bar u},
			\\
			\label{eq:special_vk_for_j3:5}
			\normLtwotime{\bar u(x) + t_k v_k(x)} - \normLtwotime{\bar u(x)} &= t_k \normLtwotime{v_k(x)}
			\quad\text{for } x \in \Omega_{\bar u}^0,
		\end{align}
	\end{subequations}
	where $K_k(x) := \normLtwotime{\bar u(x) + t_k v_k(x)} + \normLtwotime{\bar u(x)}$.
\end{lemma}
\begin{proof}
	For \eqref{eq:special_vk_for_j3:1}, we argue like in \Cref{lem:special_vk_for_j1}.
	We have pointwise convergence
	$v_k \to v$ and $v \in \LtwoOmegaT$ dominates $v_k$.
	The dominated convergence theorem yields the claim.
	
	Next we address \eqref{eq:special_vk_for_j3:3}.
	The case $\Omega_{\bar u} \times (0,T)$ can be handled as in the proof of \cref{lem:special_vk_for_j1}.
	The only interesting case is $(x,t) \in \Omega_{\bar u}^0 \times (0,T)$ if $v_k(x,t) \ne 0$.
	In this case we have
	$\bar u(x,t) = 0$, $v_k(x,t) = v(x,t)$ and $\normLtwotime{v(x)} \le \frac 1 {\sqrt{t_k}}$.
	With \eqref{subdifferential j3 charact_1}, we get
	\begin{align*}
		\abs{\bar u(x,t) + t_k v_k(x,t)}
		=
		t_k \abs{v(x,t)}
		=
		t_k \abs{\lambda_{\bar u}(x,t)} \normLtwotime{v(x)}
		\le
		\sqrt{t_k} \norm{\lambda_{\bar u}}_{L^\infty(\Omega_{\bar u}^0 \times (0,T))}.
	\end{align*}
	By combining \eqref{necessary first order j 0}
	with \eqref{eq:normal_cone},
	we get $\lambda_{\bar u} = -F'(\bar u) / \mu$
	on the set $\Omega_{\bar u}^0$.
	Together with $F'(\bar u) \in \LinftyOmegaT$,
	we get
		$\abs{\bar u(x,t) + t_k v_k(x,t)} \le C \sqrt{t_k}$
		for some constant $C > 0$.
		As $\bar u(x,t) = 0$, we have $\alpha \le 0 \le \beta$.
		If on the one hand $\alpha < 0 < \beta$, we have
		$\alpha \le \bar u(x,t) + t_k v_k(x,t) \le \beta$
		for $k$ large enough.
		If on the other hand $\alpha = 0 < \beta$ holds,
		then $0 \le \bar u(x,t) + t_k v_k(x,t)$ follows from $v \in \criticalcone{\bar u}$.
		The upper bound $\bar u(x,t) + t_k v_k(x,t) \le \beta$ holds for $k$ large enough as in the other case.
		Finally, the case $\alpha < 0 = \beta$ is similar.
	This verifies \eqref{eq:special_vk_for_j3:3}
	and we also get
	$v_k \in \tangentcone{\bar u}$.

	In order to obtain \eqref{eq:special_vk_for_j3:2},
	we use \cref{Lemma bei j1 ueber F + mu j' = 0}
	in combination with \cref{lem:equality_subdif_directional_3}.
	As in the proof of \cref{lem:special_vk_for_j1},
	we analogously get that \eqref{Lemma j zu Zeigendes:1}
	is valid with $v$ replaced by $v_k$.
	From \cref{Lemma bei j1 ueber F + mu j' = 0},
	we get that \eqref{subdifferential j3 charact_1}
	holds.
	The special definition of $v_k$ on $\Omega_{\bar u}^0 \times (0,T)$
	ensures that \eqref{subdifferential j3 charact_1}
	is also satisfied if we replace $v$ by $v_k$.
	Thus,
	$j_3'(\bar u; v_k) = \dual{\lambda_{\bar u}}{v_k}$
	and \cref{Lemma bei j1 ueber F + mu j' = 0}
	gives \eqref{eq:special_vk_for_j3:2}.
	
	The identities
	\eqref{eq:special_vk_for_j3:4}
	and
	\eqref{eq:special_vk_for_j3:5}
	hold as $L^2(0,T)$ is a Hilbert space.
\end{proof}

Finally,
the next lemma provides the convergence of some integrals.
\begin{lemma}\label{lem:j3pp_special_convergence_properties}
	We assume $-F'(\bar u) \in \partial G(\bar u) \cap \LinftyOmegaT$.
	Let $v \in \criticalcone{\bar u}$ be given
	such that
	\begin{equation}
		\label{eq:integral_squared_L2Norm_of_v__divided_by_L2Norm_of_u_is_finite}
		\int_{\Omega_{\bar u}} {\frac{\normLtwotime{v(x)}^2}{\normLtwotime{\bar u(x)}}} \dx < \infty
		.
	\end{equation}
	For a given sequence $\seq{t_k} \subset \R^+$ with $t_k \searrow 0$
	we consider the sequence
	$\seq{v_k} \subset \LtwoOmegaT$ as defined in \Cref{lem:special_vk_for_j3}.
	We further denote
	$K_k(x) := \normLtwotime{\bar u(x) + t_k v_k(x)} + \normLtwotime{\bar u(x)}$.
	Then, it holds
	\begin{subequations}
		\label{eq:j3pp_special_convergence_properties}
		\begin{align}\label{eq:j3pp_special_convergence_properties_1}
			\int_{\Omega_{\bar u}} {\frac{\int_0^T {v_k^2} \dt}{K_k(x)}} \dx
			&\to
			\int_{\Omega_{\bar u}} {\frac{\int_0^T {v^2} \dt}{2 \normLtwotime{\bar u(x)}}} \dx,
			\\
			\label{eq:j3pp_special_convergence_properties_2}
			\int_{\Omega_{\bar u}} {
			\frac{\parens*{ \int_0^T {\bar u v_k} \dt }^2}{K_k(x)^2 \normLtwotime{\bar u(x)}}} \dx
			&\to
			\int_{\Omega_{\bar u}} {
			\frac{\parens*{ \int_0^T {\bar u v} \dt }^2}{4 \normLtwotime{\bar u(x)}^3}} \dx,
			\\
			\label{eq:j3pp_special_convergence_properties_3}
			t_k \int_{\Omega_{\bar u}} {
				\frac{\int_0^T {v_k^2} \dt \int_0^T {\bar u v_k} \dt}{K_k(x)^2 \normLtwotime{\bar u(x)}}
			} \dx & \to 0.
		\end{align}
	\end{subequations}
\end{lemma}
\begin{proof}
	In the proof of \cref{lem:special_vk_for_j3},
	we have seen that
	$v_k \to v$ pointwise almost everywhere.
	Let us denote by $N \subset \OmegaT$ the null set
	on which the sequence does not converge.
	Then, we know that for almost all $x \in \Omega$,
	the set $\set{t \in (0,T) \given v_k(x,t) \not\to v(x,t) }$
	is measurable and also a null set.
	Together with $\abs{v_k} \le \abs{v}$ pointwise a.e.\ and
	$\normLtwotime{v(x)} < \infty$ for a.a.\ $x \in \Omega$
	we
	get
	$\normLtwotime{v_k(x) - v(x)} \to 0$ for a.a.\ $x \in \Omega$
	from the dominated convergence theorem.
	This shows that the integrands in \eqref{eq:j3pp_special_convergence_properties}
	converge pointwise a.e.\ on $\Omega_{\bar u}$.

	In order to apply the dominated convergence theorem,
	we only need integrable bounds.
	These can be easily obtained with $\abs{v_k} \le \abs{v}$,
	the estimates
	$K_k(x) \ge t_k \normLtwotime{v_k(x)}$,
	$K_k(x) \ge \normLtwotime{\bar u(x)}$
	and
	\eqref{eq:integral_squared_L2Norm_of_v__divided_by_L2Norm_of_u_is_finite}.
\end{proof}

\begin{theorem}\label{G3''}
	We assume $-F'(\bar u) \in \partial G_3(\bar u) \cap \LinftyOmegaT$.
	For all $v \in \criticalcone{\bar u}$ we have
	\begin{equation*}
		G_3''(\bar u, -F'(\bar u); v) =
		\begin{cases}
			\int_{\Omega_{\bar u}}
			{\frac \mu {\normLtwotime{\bar u(x)}}
				\bracks*{
					\int_0^T {v^2} \dt - \parens*{ \int_0^T {\frac{\bar u v}{\normLtwotime{\bar u(x)}}} \dt }^2
			}}
			\dx,
			&\bar u \ne 0,
			\\
			0,
			&\bar u = 0.
		\end{cases}
	\end{equation*}
	Moreover, $G_3$ is strongly twice epi-differentiable at $\bar u$ for $-F'(\bar u)$.
\end{theorem}
Note that the value $G''(\bar u, -F'(\bar u); v) = \infty$ is possible for
$\bar u \ne 0$ and $v \in \criticalcone{\bar u}$.
\begin{proof}
	We first consider the case $\bar u \ne 0$.
	We are going to use \Cref{lem:density_lemma} with
	\begin{align*}
		Q(v) :=
		\delta_{\criticalcone{\bar u}}(v)
		+
			\int_{\Omega_{\bar u}}
			{\frac \mu {\normLtwotime{\bar u(x)}}
				\bracks*{
					\int_0^T {v^2} \dt - \parens*{ \int_0^T {\frac{\bar u v}{\normLtwotime{\bar u(x)}}} \dt }^2
			}}
			\dx
	\end{align*}
	and
	$V = \set{v \in \criticalcone{\bar u} \given
		\eqref{eq:integral_squared_L2Norm_of_v__divided_by_L2Norm_of_u_is_finite} \text{ holds.}}$.
	We have to check the assumptions of \cref{lem:density_lemma}.
	
	Step 1, \itemref{lem:density_lemma:1}:
	For $v \not\in \criticalcone{\bar u}$,
	we have 
	$G_3''(\bar u, -F'(\bar u); v) = \infty = Q(v)$,
	see \cref{lem:Gpp_directional_derivative}.
	Let $v \in \criticalcone{\bar u}$ be arbitrary and consider
	sequences $\seq{t_k} \subset \R^+$ and $\seq{v_k} \subset \LtwoOmegaT$
	with $t_k \searrow 0$ and $v_k \weakly v$.
	We use $-F'(\bar u) \in \partial G_3(\bar u)$
	to get
	$0 \le \dual{F'(\bar u)}{v_k} + G_3'(\bar u; v_k) = \dual{F'(\bar u)}{v_k} + \indicatorofUad'(\bar u;v_k) + \mu j_3'(\bar u; v_k)$.
	Consequently,
	\begin{align*}
		\MoveEqLeft
		\liminf_{k \to \infty} \frac{G_3(\bar u + t_k v_k) - G_3(\bar u) - t_k \dual{-F'(\bar u)}{v_k}}{t_k ^2 / 2}
		\\
		&
		=
		\liminf_{k \to \infty} \frac 2 {t_k^2} \parens*{
			\mu \bracks*{ j_3(\bar u + t_k v_k) - j_3(\bar u) }
			+
			\indicatorofUad(\bar u + t_k v_k)
		+ t_k \dual{F'(\bar u)}{v_k} }
		\\
		&
		\ge
		\liminf_{k \to \infty} \frac {2} {t_k^2} \parens[\bigg]{
			\mu \bracks*{ j_3(\bar u + t_k v_k) - j_3(\bar u) }
		+ t_k \braces*{ \indicatorofUad'(\bar u;v_k) + \dual{F'(\bar u)}{v_k} }
		}.
		\\
		&
		\ge
		\liminf_{k \to \infty} \frac {2 \mu} {t_k^2} \parens*{ j_3(\bar u + t_k v_k) - j_3(\bar u) - t_k j_3'(\bar u; v_k) }
		\\
		&
		\ge
		\int_{\Omega_{\bar u}} 
		{\frac \mu {\normLtwotime{\bar u(x)}} \braces*{ \int_0^T {v(x,t)^2} \dt - \frac {\parens*{ \int_0^T {\bar u(x,t) v(x,t)} \dt }^2} {\normLtwotime{\bar u(x)}^2} }}
		\dx
		=
		Q(v),
	\end{align*}
	where we used \eqref{eq:j3_first_part_of_proof:estimate_for_j3pp-term}
	in the last step.
	Taking the infimum with respect to the sequences $\seq{t_k}$, $\seq{v_k}$
	yields $G_3''(\bar u, -F'(\bar u); v) \ge Q(v)$ for all $v \in \LtwoOmegaT$.

	Step 2, \itemref{lem:density_lemma:2}:
	We consider arbitrary $v \in V$ and $\seq{t_k} \subset \R^+$ with $t_k \searrow 0$. Let
	$\seq{v_k} \subset \LtwoOmegaT$ be defined as in \cref{lem:special_vk_for_j3}.
	We have
	\begingroup\allowdisplaybreaks
	\begin{align*}
		&\lim_{k \to \infty} \frac{G_3(\bar u + t_k v_k) - G_3(\bar u) - t_k \dual{-F'(\bar u)}{v_k}}{t_k ^2 / 2}
		\\
		&= \lim_{k \to \infty} \frac 2 {t_k^2} \parens*{ \mu \bracks*{ j_3(\bar u + t_k v_k) - j_3(\bar u) } - t_k \mu j_3'(\bar u; v_k) }
		&&\mathllap{\text{(by \eqref{eq:special_vk_for_j3:2}, \eqref{eq:special_vk_for_j3:3})}}
		\\
		&= \lim_{k \to \infty} \frac {2 \mu} {t_k^2}
		\int_{\Omega_{\bar u}}
		\normLtwotime{\bar u + t_k v_k} - \normLtwotime{\bar u}
		- t_k
		\frac {\int_0^T {\bar u v_k} \dt} {\normLtwotime{\bar u}}
		\dx
		&&\mathllap{\text{(by \eqref{eq:directional_derivative_j3}, \eqref{eq:special_vk_for_j3:5})}}
		\\
		&= \lim_{k \to \infty} \frac {2 \mu} {t_k^2}
		\int_{\Omega_{\bar u}} {\frac{\int_0^T {2 t_k \bar u v_k + t_k^2 v_k^2} \dt}{K_k} - \frac {t_k} {\normLtwotime{\bar u}} \int_0^T {\bar u v_k} \dt} \dx
		&&\mathllap{\text{(by \eqref{eq:special_vk_for_j3:4})}}
		\\
		&= \lim_{k \to \infty} 2 \mu
		\int_{\Omega_{\bar u}} {\frac{\int_0^T {v_k^2} \dt}{K_k} +
			\frac 1 {t_k} \parens*{
				\frac{2}{K_k} - \frac {1} {\normLtwotime{\bar u}}
			}
		\int_0^T {\bar u v_k} \dt} \dx
		\\
		&= \lim_{k \to \infty} 2 \mu
		\int_{\Omega_{\bar u}} {\frac{\int_0^T {v_k^2} \dt}{K_k} +
			\frac 1 {t_k} \parens*{
				\frac{\normLtwotime{\bar u} - \normLtwotime{\bar u + t_k v_k}}{K_k \normLtwotime{\bar u}}
			}
		\int_0^T {\bar u v_k} \dt} \dx
		\qquad\qquad
		\\
		&= \lim_{k \to \infty} 2 \mu
		\int_{\Omega_{\bar u}} {\frac{\int_0^T {v_k^2} \dt}{K_k} -
			\parens*{
				\frac{\int_0^T {2 \bar u v_k + t_k v_k^2} \dt}{K_k^2 \normLtwotime{\bar u}}
			}
		\int_0^T {\bar u v_k} \dt} \dx
		&&\mathllap{\text{(by \eqref{eq:special_vk_for_j3:4})}}
		\\
		&= \lim_{k \to \infty} 2 \mu
		\int_{\Omega_{\bar u}} {\frac{\int_0^T {v_k^2} \dt}{K_k}
			-
			\frac{2 \parens*{ \int_0^T {\bar u v_k} \dt }^2}{K_k^2 \normLtwotime{\bar u}}
			-
			t_k \frac{\int_0^T {v_k^2} \dt \int_0^T {\bar u v_k} \dt}{K_k^2 \normLtwotime{\bar u}}
		} \dx
		\\
		&=
		\int_{\Omega_{\bar u}} {\frac \mu {\normLtwotime{\bar u}} \bracks*{ \int_0^T {v^2} \dt - \parens*{ \int_0^T {\frac{\bar u v}{\normLtwotime{\bar u}}} \dt }^2 }} \dx.
		&&\mathllap{\text{(by \eqref{eq:j3pp_special_convergence_properties})}}
	\end{align*}
	\endgroup

	Step 3, \itemref{lem:density_lemma:3}:
	For arbitrary $v \in \criticalcone{\bar u}$ we define the sequence $\seq{v^l} \subset \LtwoOmegaT$ via
	\begin{equation*}
		v^l(x,t) :=
		\begin{cases}
		0      &\text{if } 0 < \normLtwotime{\bar u(x)} < \frac 1 l, \\
		v(x,t) &\text{else}.
		\end{cases}
	\end{equation*}
	It is clear that $v^l \to v$ in $\LtwoOmegaT$
	and $v^l \in \tangentcone{\bar u}$.
	We are going to use \cref{lem:equality_subdif_directional_1,lem:equality_subdif_directional_3}
	to check $v^l \in \criticalcone{\bar u}$.
	The property \eqref{Lemma j zu Zeigendes:1} holds for $v^l$ by construction.
	\cref{lem:equality_subdif_directional_3} yields \eqref{subdifferential j3 charact_1} for $v$.
	Now, it is straightforward to check that
	\eqref{subdifferential j3 charact_1} also holds for $v$ replaced by $v^l$.
	Therefore, \cref{lem:equality_subdif_directional_3} yields
	$\dual{\lambda_{\bar u}}{v^l} = j_3'(\bar u; v^l)$.
	\Cref{Lemma bei j1 ueber F + mu j' = 0} yields $v^l \in \criticalcone{\bar u}$ for all $l \in \N$.
	
	Due to
	\begin{equation*}
		\int_{\Omega_{\bar u}} {\frac {\normLtwotime{v^l(x)}^2}{\normLtwotime{\bar u(x)}}} \dx
		=
		\int_{\Omega_{1/l}} {\frac {\normLtwotime{v(x)}^2}{\normLtwotime{\bar u(x)}}} \dx
		\le
		l \Ltwonorm{v}^2
		<
		\infty
		,
	\end{equation*}
	we have $v^l \in V$.
	
	From \eqref{eq:integrand_of_j3_nonnegative}
	we get $Q(v) \ge Q(v^l)$
	and therefore $Q(v) \ge \liminf_{l \to \infty} Q(v^l)$.
	
	Now, we are in position to apply
	\cref{lem:density_lemma}
	and this yields the claim
	in case $\bar u \ne 0$.

	Finally, it remains to consider the case $\bar u = 0$.
	\Cref{lem:Gpp_convex} yields $G_3''(0,-F'(0);v) \ge 0$.
	We consider an arbitrary sequence $\seq{t_k} \subset \R^+$ with $t_k \searrow 0$
	and choose $\seq{v_k} \subset \LtwoOmegaT$ as in \cref{lem:special_vk_for_j3}.
	We get
	\begin{align*}
		\lim_{k \to \infty} \frac{G_3(0 + t_k v_k) - G_3(0) - t_k \dual{-F'(0)}{v_k}}{t_k ^2 / 2}
		&=
		\lim_{k \to \infty} \frac {\mu t_k j_3(v_k) + t_k \dual{F'(0)}{v_k}} {t_k^2 / 2} \\
		&=
		\lim_{k \to \infty} \frac {\mu j_3'(0;v_k) + \dual{F'(0)}{v_k}} {t_k / 2} \\
		&=
		0,
	\end{align*}
	using
	\eqref{eq:special_vk_for_j3:3},
	$j_3(v_k) = j_3'(0;v_k)$ (see \eqref{eq:directional_derivative_j3}) and
	\eqref{eq:special_vk_for_j3:2}.
	This finishes the proof.
\end{proof}

\section{Application to a parabolic control problem}
\label{sec:application}

In this section,
we apply the findings from \cref{sec:subderivatives_sparsity}
to the optimal control problem
\begin{equation}
	\tag{OCP}
	\label{eq:optimal_control_problem_for_special_F}
	\begin{aligned}
		\text{Minimize} \quad &J(u) = F(u) + \mu j(u),
		\\
		\text{w.r.t.} \quad & u \in \Uad,
	\end{aligned}
\end{equation}
where the smooth part $F$
is given by
\begin{align*}
	F(u) = \int_{\OmegaT} L(x,t,y_u(x,t)) \dxt +
	\frac \nu 2 \Ltwonorm{u}^2
\end{align*}
and $y_u \in W(0,T) := \set{y \in L^2(0,T; H_0^1(\Omega)) \given \partial_t y \in L^2(0,T;H^{-1}(\Omega))}$
is the (weak) solution of the state equation
\begin{equation}\label{DGL yu}
			\partial_t y_u + A y_u + a(\cdot, y_u) = u \text{ in } \OmegaT,
			\qquad
			y_u = 0 \text{ on } \Sigma_T, \qquad
			y_u(\cdot,0) = y_0 \text{ in } \Omega
\end{equation}
As in \cref{sec:subderivatives_sparsity},
$\OmegaT := \Omega \times (0,T)$,
where
$\Omega \subset \R^d$ is assumed to be non-empty, open, and bounded,
and $T > 0$.
Moreover, $\Sigma_T := \Gamma \times (0, T)$.
The nonsmooth part $j$ is one of the functionals in \eqref{eq:js}
and we define $G$ as in \eqref{eq:G},
i.e., $G := \delta_{\Uad} + \mu j$.
Note that $\dom(G) = \Uad$.
We further assume $\mu > 0$, $\nu \ge 0$
and the bounds $\alpha,\beta \in \R$
satisfy $\alpha < \beta$.

The above control problem has been analyzed
in \cite{CasasHerzogWachsmuth2015:1}.
In order to compare our results,
we rely on the same standing assumptions,
which are assumed to hold throughout this section.
\begin{assumption}
	\label{asm:assumptions_for_control_problem}
	We assume that $A$, $a$, $L$
	together with exponents
	$\hat p, \hat q \in [2, \infty]$
	satisfy
	\cite[Assumptions~1--3]{CasasHerzogWachsmuth2015:1}.
	In particular, $A$ is an elliptic differential operator
	and the Nemytskii operators $a$ and $L$ satisfy
	the usual continuity, differentiability and growth conditions
	depending on $\hat p$ and $\hat q$.
\end{assumption}
From this assumption,
we get the following two differentiability results.
\begin{lemma}[\texorpdfstring{\cite[Theorem~2.1]{CasasHerzogWachsmuth2015:1}}{[Theorem 2.1, Casas et al., 2017]}]
	\label{Lemma DGL zv und zv1v2}
	For all $u \in L^{\hat p}(0,T;L^{\hat q}(\Omega))$ the equation
	\eqref{DGL yu}
	has a unique solution
	$y_{u} \in W(0,T) \cap L^\infty(\OmegaT)$.
	Moreover, the solution mapping
	$H: L^{\hat p}(0,T;L^{\hat q}(\Omega)) \to W(0,T) \cap L^\infty(\OmegaT)$, defined by
	$H(u) := y_{u}$, is of class $C^2$.
	For all elements $u, v \in L^{\hat p}(0,T;L^{\hat q}(\Omega))$,
	the function $z_v = H'(u)v$
	is the solutions of the problem
	\begin{equation}\label{DGL zv}
		\frac{\partial z}{\partial t} + Az + \frac{\partial a}{\partial y}(\cdot,y_{u}) z = v \text{ in } \OmegaT,
		\qquad
		z = 0 \text{ on } \Sigma_T,
		\qquad
		z(\cdot,0) = 0 \text{ in } \Omega,
	\end{equation}
	respectively.
\end{lemma}

\begin{lemma}[\texorpdfstring{\cite[Theorem~2.3]{CasasHerzogWachsmuth2015:1}}{[Theorem 2.3, Casas et al., 2017]}]
	\label{Form von F F' und F''}
	The map
	$F: L^{\hat p}(0,T;L^{\hat q}(\Omega)) \to \R$
	is of class $C^2$. Moreover, for all
	$u, v, v_1, v_2 \in L^{\hat p}(0,T;L^{\hat q}(\Omega))$
	we have
	\begin{subequations}
		\label{F' und F'' Form}
		\begin{align}
			\label{F' Form}
			F'(u) v
			&=
			\int_\OmegaT {(\varphi_{u} + \nu u) v} \dxt
			\\
			\label{F'' Form}
			F''(u)(v_1, v_2) &= \int_\OmegaT
			\braces*{
				\parens*{
					\frac{\partial^2 L}{\partial y^2}(x,t,y_{u})
					- \varphi_{u} \frac{\partial^2 a}{\partial y^2}(x,t,y_{u})
				}  z_{v_1} z_{v_2}
				+ \nu v_1 v_2
			}
			\dxt,
		\end{align}
	\end{subequations}
	where $z_{v_i} = H'(u) v_i, i = 1, 2$, and $\varphi_{u} \in W(0,T) \cap L^\infty(\OmegaT)$
	is the solution of
	\begin{equation*}
		-\frac{\partial \varphi}{\partial t} + A^* \varphi  + \frac{\partial a}{\partial y}(\cdot,y_{\bar u}) \varphi = \frac{\partial L}{\partial y}(\cdot,y_{\bar u}),
		\qquad
		\varphi = 0 \text{ on } \Sigma_T,
		\qquad
		\varphi(\cdot, T) = 0 \text{ in } \Omega,
	\end{equation*}
	where $A^*$ is the adjoint operator of $A$.
\end{lemma}

In view of \cref{G3''},
we note that
$\varphi_{u}, u \in L^\infty(\OmegaT)$
implies $F'(u) \in L^\infty(\OmegaT)$.

For later reference,
we state the following very important estimate
and the compactness of the mapping $v \mapsto z_v$.

\begin{lemma}
	\label{lem:vk_weakly_towards_v => zvk_strongly_towards_zv}
	Let $\bar u \in \Uad$ be given.
	Then, there exists $C_Z > 0$ satisfying
	\begin{equation}\label{Normabschaetzung Anwendungsfall Lemma}
		\Ltwonorm{z_{v}} \le C_Z \Ltwonorm{v}
		\qquad
		\forall v \in \LtwoOmegaT
		.
	\end{equation}
	Additionally, if $v_k \weakly v$ in $\LtwoOmegaT$ holds, then $z_{v_k} \to z_v$ in $\LtwoOmegaT$.
\end{lemma}

The norm estimate \eqref{Normabschaetzung Anwendungsfall Lemma}
follows from standard parabolic estimates
and the compactness is a consequence
of the celebrated Aubin--Lions lemma, see, e.g., \cite{Simon1986}.

Finally,
we cite the next continuity result
for the second derivative of $F$.
\begin{lemma}
	[\texorpdfstring{\cite[Lemma~5.2]{CasasHerzogWachsmuth2015:1}}{[Lemma 5.2, Casas et al., 2017]}]
	\label{Lemma zu Voraussetzungen im Anwendungsfall}
	Let $\bar u \in \Uad$ be given.
	For all $\rho > 0$ there exists $\varepsilon > 0$ such that
	\begin{equation*}
		\abs{F''(u) v^2 - F''(\bar u) v^2} \le \rho \Ltwonorm{z_v}^2
		\qquad
		\forall v \in \LtwoOmegaT, u \in \Uad, \Ltwonorm{u - \bar u} \le \varepsilon
		,
	\end{equation*}
	where $z_v = H'(\bar u)v$.
\end{lemma}

Now, we are in position
to check that \cref{asm:standing_assumption}
is satisfied.

\begin{theorem}
	\label{Anwendungsfall Voraussetzungen erfuellt}
	The functional $F$ maps $\dom(G) = \Uad$ to $\R$.
	Let $\bar u \in \Uad$ be fixed.
	The (bi)linear functionals $F'(\bar u)$ and $F''(\bar u)$
	defined on 
	$L^{\hat p}(0,T;L^{\hat q}(\Omega))$
	can be extended to continuous (bi)linear functionals on $\LtwoOmegaT$.
	Then, $F$ satisfies the assumptions in \Cref{asm:standing_assumption}
	with $X = \LtwoOmegaT$ and $\bar x = \bar u \in \Uad$.
\end{theorem}
\begin{proof}
	Since $X = \LtwoOmegaT$
	is a separable Hilbert space,
	\itemref{asm:standing_assumption:1}
	holds.

	It remains to check
	\itemref{asm:standing_assumption:3}.
	We already mentioned
	that $\varphi_{\bar u}, \bar u \in \LinftyOmegaT$,
	thus
	$F'(\bar u) \in \LinftyOmegaT \subset \LtwoOmegaT$.
	Next, we check that
	the bilinear form $F''(\bar u)$
	can be extended to $\LtwoOmegaT$.
	Due to the assumptions made on $L$ and $a$,
	one can show that the term in parentheses
	in \eqref{F'' Form}
	belongs to $\LinftyOmegaT$.
	Thus,
	together with \cref{lem:vk_weakly_towards_v => zvk_strongly_towards_zv},
	we get
	\begin{equation}\label{Abschaetzung F'' Anwendungsfall}
		\begin{aligned}
			\abs{F''(\bar u)(v_1, v_2)}
			&\le
			C \Ltwonorm{z_{v_1}} \Ltwonorm{z_{v_2}} + \nu \Ltwonorm{v_1} \Ltwonorm{v_2} \\
			&\le
			\parens[\big]{ C C_Z^2 + \nu } \Ltwonorm{v_1} \Ltwonorm{v_2}.
		\end{aligned}
	\end{equation}
	Together with the density
	of
	$L^{\hat p}(0,T;L^{\hat q}(\Omega))$
	in
	$\LtwoOmegaT$,
	we can extend $F''(\bar u)$ continuously to $\LtwoOmegaT$.

	We still have to check \eqref{eq:hadamard_taylor_expansion}.
	Let $\seq{t_k} \subset \R^+$ and $\seq{v_k} \subset \LtwoOmegaT$
	with $t_k \searrow 0$, $v_k \weakly v \in \LtwoOmegaT$
	and
	$\bar{u} + t_k v_k \in \dom (G) = \Uad$ be given.
	For an arbitrary $\rho > 0$,
	we utilize
	\cref{Lemma zu Voraussetzungen im Anwendungsfall,lem:vk_weakly_towards_v => zvk_strongly_towards_zv}
	(together with $t_k \searrow 0$ and the boundedness of $\seq{v_k}$ in $\LtwoOmegaT$)
	to get
	\begin{equation*}
		\abs{F''(\bar u + \theta_k t_k v_k) v_k^2 - F''(\bar u) v_k ^2}
		\le
		\rho
	\end{equation*}
	for all $\theta_k \in [0,1]$ and all $k$ large enough (depending on $\rho$).
	Next, we use a second-order Taylor expansion
	and obtain intermediate points $\seq{\theta_k} \subset [0,1]$
	such that
	\begin{align*}
		\abs*{\frac{F(\bar u + t_k v_k) - F(\bar u) - t_k F'(\bar u)v_k - \frac12 t_k ^2 F''(\bar u) v_k ^2}{t_k ^2}}
		&=
		\frac 1 2
		\abs{F''(\bar u + \theta_k t_k v_k) v_k^2 - F''(\bar u) v_k ^2} \\
		&\le
		\rho
	\end{align*}
	for all $k$ large enough (depending on $\rho$).
	Since $\rho > 0$ was arbitrary,
	this shows \eqref{eq:hadamard_taylor_expansion}.
\end{proof}

Moreover,
the following holds for the nonsmooth part of the objective.
As in \cref{sec:subderivatives_sparsity},
the functional $G$
can represent any of the functionals $G_i$, $i = 1,2,3$.
\begin{lemma}
	\label{lem:second_order_nonsmooth}
	Let $\bar u$ be given such that $-F'(\bar u) \in \partial G(\bar u)$.
	In case $j = j_2$, we additionally assume $\bar u \ne 0$.
	Then,
	the functional $G$ is strongly twice epidifferentiable
	at $\bar u$ w.r.t.\ $-F'(\bar u)$.
	Moreover, $G''(\bar u, -F'(\bar u); v) = \infty$ for all $v \in \LtwoOmegaT \setminus \criticalcone{\bar u}$.
\end{lemma}
\begin{proof}
	The strong twice epidifferentiability
	follows from \cref{G1'',G2'',G3''}.
	In case $G = G_3$ we also need $-F'(\bar u) \in \LinftyOmegaT$,
	which was provided after \cref{Form von F F' und F''}.
	The final assertion is \eqref{eq:heile_welt_on_KK}.
\end{proof}

Now we can prove the second-order necessary conditions.
Recall that expressions for
$G''(\bar u, -F'(\bar u); \cdot)$
where given in
\cref{G1'',G2'',G3''}.
\begin{theorem}[Second-Order Necessary Conditions]
	\label{thm:SNC_of_OCP}
	Let $\nu \ge 0$ be given and
	let $\bar u \in \LtwoOmegaT$ be a local minimizer of \eqref{eq:optimal_control_problem_for_special_F}.
	In case $j = j_2$, we additionally assume $\bar u \ne 0$.
	Then,
	$-F'(\bar u) \in \partial G(\bar u)$
	and
	$F''(\bar u) v^2 + G''(\bar u,-F'(\bar u);v) \ge 0$ holds for all $v \in \criticalcone{\bar u}$.
\end{theorem}
\begin{proof}
	We apply \cref{thm:SNC} with the setting
	$X = \LtwoOmegaT$,
	$\bar x := \bar u$,
	$F := F$,
	$G := \indicatorofUad + \mu j_i$, for some $i = 1,2,3$,
	$c := 0$.
	Hence, \eqref{eq:second_order_growth} holds.
	The first-order condition \cref{thm:FONC}
	yields
	$-F'(\bar u) \in \partial G(\bar u)$.
	Additionally, \cref{thm:SNC}~\ref{assumption-mrc} is satisfied,
	see \cref{lem:second_order_nonsmooth}.
\end{proof}
Let us compare this result with
\cite[Theorem~4.3]{CasasHerzogWachsmuth2015:1}.
We have already seen that the critical cones
in both results coincide, see \cref{lem:critical_cones}.
Further,
it can be checked that
the expressions
for $G''(\bar u, -F'(\bar u); v)$
given in
\cref{G1'',G2'',G3''}
coincide with the corresponding expressions for
$j_i''(\bar u; v^2)$
given in \cite[(4.5)--(4.7)]{CasasHerzogWachsmuth2015:1}.
Hence, both results coincide.
Note that \cite[Theorem~4.3]{CasasHerzogWachsmuth2015:1}
also addresses the case $\bar u = 0$ if $j = j_2$,
but the proof is flawed.
Indeed,
it is claimed (without any justification) that for arbitrary $v \in \criticalcone{\bar u}$
we have $v_k := \projection{[-k,k]}{v} \in \criticalcone{\bar u}$,
but \cref{ex:critical_cone_j2}
shows that this might fail.

In order to apply the second-order sufficient assumptions,
we need an additional lemma.
\begin{lemma}\label{lem:Fpp(u)_is_seq_weakstar_semicts}
	In case $\nu > 0$, the mapping
	\begin{equation*}
		\LtwoOmegaT \ni v \mapsto F''(\bar u) {v}^2 =
		\int_\OmegaT
		{
			\parens*{
			\frac{\partial^2 L}{\partial y^2}(x,t,y_{\bar u})
			- \varphi_{\bar u} \frac{\partial^2 a}{\partial y^2}(x,t,y_{\bar u})
		}  z_{v}^2
			}
		\dxt
		+ \nu \Ltwonorm{v}^2
	\end{equation*}
	is a Legendre form.
\end{lemma}
\begin{proof}
	Due to the compactness result
	\cref{lem:vk_weakly_towards_v => zvk_strongly_towards_zv},
	the first addend in $F''(\bar u) v^2$ is
	sequentially weakly continuous.
	Applying
	\cite[Proposition~3.76]{BonnansShapiro2000}
	yields the claim.
\end{proof}
In case $\nu = 0$, the map $v \mapsto F''(\bar u) v^2$ is sequentially weakly continuous
and, thus, not a Legendre form.

\begin{theorem}[Second-Order Sufficient Condition]
	We assume $\nu > 0$.
	Further suppose that $\bar u \in \Uad$
	satisfies $-F'(\bar u) \in \partial G(\bar u)$
	and
	\begin{equation}\label{Hinr. 2. Ordnung Anwendungsfall Voraussetzung J'' > 0}
		F''(\bar u)v^2 + G''(\bar u,-F'(\bar u);v) > 0
		\qquad
		\forall v \in \criticalcone{\bar u} \setminus \set{ 0 }
		.
	\end{equation}
	Then, there exist $\varepsilon, \delta > 0$ such that
	\begin{equation}
		J(u) \ge J(\bar u) + \frac \delta 4 \Ltwonorm{u - \bar u}^2
		\qquad
		\forall
		u \in \Uad, \Ltwonorm{u - \bar u} \le \varepsilon.
	\end{equation}
\end{theorem}
\begin{proof}
	We will show the requirements for the application of \Cref{thm:SSC_wo}.
	\Cref{lem:Fpp(u)_is_seq_weakstar_semicts} yields the sequential weak lower semicontinuity of $v \mapsto F''(\bar u) {v}^2$.
	Condition \eqref{eq:SSC_wo} follows from
	\eqref{Hinr. 2. Ordnung Anwendungsfall Voraussetzung J'' > 0}
	and $-F'(\bar u) \in \partial G(\bar u)$,
	cf.\ the proof of \cref{cor:2nd_order_traditional_noncvx}.
	\Cref{lem:Fpp(u)_is_seq_weakstar_semicts}
	in combination with
	\cref{lem:Legendre_gives_NDC}
	shows that 
	\eqref{eq:NDC}
	holds.
	Now, the claim follows from \cref{thm:SSC_wo}.
\end{proof}
Let us compare this result with the second-order results
in \cite[Section~5]{CasasHerzogWachsmuth2015:1}
in the case $\nu > 0$.
For the functional $j_1$
we obtain an identical result.

For the functional $j_2$,
we get the same result in case $\bar u \ne 0$.
Note that \cref{thm:SSC_wo} is still applicable
in case $\bar u = 0$,
although we do not know the precise values of $G_2''(\bar u,-F'(\bar u);\cdot)$.
Our sufficient condition reads
\begin{equation*}
	F''(\bar u) v^2 + G''(\bar u, -F'(\bar u); v) > 0
	\qquad \forall v \in \criticalcone{\bar u} \setminus \set{0}
\end{equation*}
and due to $G''(\bar u, -F'(\bar u); v) \ge 0$,
this condition is weaker than the sufficient condition
\begin{equation*}
	F''(\bar u) v^2 > 0
	\qquad \forall v \in \criticalcone{\bar u} \setminus \set{0}
\end{equation*}
given in
\cite[Theorem~5.8]{CasasHerzogWachsmuth2015:1},
see also \cite[(4.6)]{CasasHerzogWachsmuth2015:1}.

For the functional $j_3$,
\cite[Theorem~5.12]{CasasHerzogWachsmuth2015:1}
shows the quadratic growth
only in an $L^\infty(\Omega; L^2(0,T))$-ball,
whereas our result
implies the growth in the larger $\LtwoOmegaT$-ball.

We also note that
\cite[Section~5]{CasasHerzogWachsmuth2015:1}
contains sufficient optimality conditions in case $\nu = 0$
which were not under investigation here.
In fact,
in this case,
the theory from \cref{sec:second-order_conditions}
cannot be applied in the space
$\LtwoOmegaT$,
since \eqref{eq:NDC}
cannot be satisfied.
For a similar problem with $\mu = 0$,
i.e., without a sparsity-inducing term,
it was shown
in \cite{ChristofWachsmuth2017:1}
that
(under a certain regularity assumption)
a second-order analysis
can be performed in a space of measures.
In \cite[Section~5.1]{WachsmuthWachsmuth2022},
this analysis was extended to $\mu > 0$
in case of the functional $j_1$.
The extension to $j_2$ and $j_3$ is subject to future work.

\bibliographystyle{jnsao}
\bibliography{references}

\end{document}